\documentclass[12pt]{amsart}
\usepackage{amsthm, amssymb}
\usepackage{amsfonts, amscd}
\usepackage{epsfig,multicol}
\theoremstyle{plain} \numberwithin{equation}{section}
\newtheorem{thm}{Theorem}[section]
\newtheorem{cor}[thm]{Corollary}
\newtheorem{lem}{Lemma}[section]
\newtheorem{prop}[thm]{Proposition}

\newtheorem{exam}{Example}[section]
\theoremstyle{definition}
\newtheorem{defn}{Definition}[section]
\newtheorem{rem}{Remark}[section]

\def\H{\text{Hom}}

\begin{document}
\title[2-torus manifolds, cobordism and small covers]
{\large \bf 2-torus manifolds, cobordism and small covers}
\author[Zhi L\"u]{Zhi L\"u }
\footnote[0]{ 
\endgraf
\endgraf
 Supported by grants from NSFC (No. 10371020 and No. 10671034). }
\address{Institute of Mathematics, School of Mathematical Sciences, Fudan University, Shanghai,
200433, P.R. China.} \email{zlu@fudan.edu.cn}
\begin{abstract}
Let ${\frak M}_n$ be the set of equivariant unoriented cobordism
classes of all $n$-dimensional 2-torus manifolds, where an
$n$-dimensional 2-torus manifold $M$ is a smooth closed manifold
of dimension $n$ with effective smooth action of a rank $n$
2-torus group $({\Bbb Z}_2)^n$. Then ${\frak M}_n$ forms an
abelian group with respect to disjoint union. This paper
determines the group structure of ${\frak M}_n$ and shows that
each class of ${\frak M}_n$ contains a small cover as its
representative in the case $n=3$.
\end{abstract}
\maketitle

\section{Introduction}

An $n$-dimensional 2-torus manifold $M$ is a smooth closed manifold
of dimension $n$ with effective smooth action of a rank $n$ 2-torus
group $({\Bbb Z}_2)^n$. Since the action is effective, the fixed
point set of the action is 0-dimensional (i.e., it is formed by
finitely many isolated points) if $M$ has a fixed point.  In this
paper, we shall study this class of geometrical objects from the
viewpoint of cobordism.

\vskip .2cm

Let ${\frak M}_n$ denote the set of equivariant unoriented
cobordism classes of all $n$-dimensional 2-torus manifolds. Then
${\frak M}_n$ forms an abelian group with respect to disjoint
union, and in particular, ${\frak M}_n$ also forms a vector space
over ${\Bbb Z}_2$. The zero element of ${\frak M}_n$ is given by a
canonical 2-torus manifold, which is the $n$-dimensional standard
sphere $S^n$ with the standard $({\Bbb Z}_2)^n$-action defined by
$$(x_0,x_1,...,x_n)\longmapsto (x_0, g_1x_1,...,g_nx_n),$$
 fixing two isolated points with same  $({\Bbb Z}_2)^n$-representation,
 where $(x_0,x_1,...,x_n)\in S^n$ and $(g_1,...,g_n)\in ({\Bbb
Z}_2)^n$. When $n=1, 2$, it is known from the work of Conner and
Floyd \cite{cf} that  ${\frak M}_1$ is trivial and ${\frak M}_2$ is
generated by  the standard $({\Bbb Z}_2)^2$-action on ${\Bbb R}P^2$.
As for $n\geq 3$, as far as the author knows, the structure of
${\frak M}_n$ is still far from well understood. One of main
objectives of this paper considers the following problem.

\vskip .2cm

{\bf Problem:} {\em To determine the group structure of ${\frak
M}_n$ when $n\geq 3$.}
 \vskip .2cm

 In 1991,  Davis and Januszkiewicz introduced and studied a kind
 of special 2-torus manifolds---small covers, each of which is
 locally isomorphic to a faithful representation of $({\Bbb
 Z}_2)^n$ on ${\Bbb R}^n$, and its orbit space is a simple convex
 polytope. This establishes a direct link between equivariant topology
 and  combinatorics. A typical example of an equivariant
  nonbounding small cover is a real projective space ${\Bbb R} P^n$ with
a standard action of $({\Bbb Z}_2)^n$.  Its orbit space is an
$n$-simplex. Another typical example of a bounding small cover is
a product of $n$ copies of a circle $S^1$ with reflection, and its
orbit space is an $n$-cube.  Thus, we see that when $n=2$, two
typical examples above can be used as representatives of two
classes in ${\frak M}_2$, respectively. This leads us to another
objective of this paper, i.e., the following conjecture.

\vskip .2cm {\bf Conjecture:} {\em Each class of ${\frak M}_n$
contains a small cover as its representative. }\vskip .2cm

Note that in non-equivariant case, the above conjecture has been
shown to be true by Bukhshtaber and Ray in \cite{br}.

\vskip .2cm

In this paper we  settle the above problem and conjecture in
3-dimensional case, see Theorems~\ref{dim} and \ref{small}.

\vskip .2cm

The paper is is organized as follows. In Section 2, we formulate the
complete equivariant cobordism invariant (i.e., prime tangent
representation set $\mathcal{N}_\beta$) of 2-torus manifolds from
Stong homomorphism, and then we study some properties of the
complete equivariant cobordism invariant $\mathcal{N}_\beta$. In
Section 3 we introduce the notion of an essential generator of
${\frak M}_n$, and show that any element of ${\frak M}_n$ is a
linear combination of essential generators. In Section 4, we review
the work of Davis and Januszkiewicz \cite{dj} and give two kinds of
3-dimensional small covers, which play a key role in the study of
${\frak M}_3$. In Section 5 we introduce the moment graphs induced
by 2-torus manifolds. The group structure of ${\frak M}_3$ is
determined completely in Section 6, and the above conjecture is
settled in the 3-dimensional case in Section 7.

\vskip .2cm

The author expresses his gratitude to Professor M. Masuda for his
valuable suggestions and comments, and especially for helpful
conversation in the argument of Proposition~\ref{bound}. The author
also expresses his gratitude to Professor R.E. Stong for his
valuable suggestions and comments.

\section{$G$-representations and Stong homomorphism}

Let $G=({\Bbb Z}_2)^n$, and let $\H(G,{\Bbb Z}_2)$ be the set of all
homomorphisms $\rho: G\longrightarrow {\Bbb Z}_2$, which consists of
$2^n$ distinct homomorphisms. One agrees to let $\rho_0$ denote the
trivial element in $\H(G,{\Bbb Z}_2)$, i.e., $\rho_0(g)=1$ for all
$g\in G$. The irreducible real $G$-representations are all
one-dimensional and correspond to all elements in $\H(G,{\Bbb
Z}_2)$. Thus, every irreducible real representation of $G$ has the
form $\lambda_\rho: G\times{\Bbb R}\longrightarrow{\Bbb R}$ with
$\lambda_\rho(g,x)=\rho(g)\cdot x$ for $\rho\in\H(G,{\Bbb Z}_2)$.

\vskip .2cm Given an element $\beta$ of ${\frak M}_n$, let $(M,
\phi)$ be a representative of $\beta$ such that $M$ has a fixed
point. Taking an isolated point $p$ in $M^G$, the $G$-representation
at $p$ can be written as
$$\tau_pM=\bigoplus_{\rho\not= \rho_0}\lambda_\rho^{q_\rho}$$
with $\sum_{\rho\not=\rho_0}q_\rho=n$. By the Borel Theorem (see
\cite{ap}) and the effectiveness of the action, if $q_\rho\not=0$,
then $q_\rho$ must be one. Thus, $\tau_pM$ is the direct sum of $n$
irreducible real $G$-representations (which are linearly
independent). The collection $\mathcal{N}_M=\{[\tau_pM] \big\vert\
p\in M^G\}$ is called the {\em tangent representation set of $(M,
\phi)$}, where $[\tau_pM]$ denotes the isomorphism class of
$\tau_pM$.

\vskip .2cm

 By $R_n(G)$ denote the vector space over ${\Bbb Z}_2$,
generated by the representation classes of dimension $n$. Then
$R_*(G)=\sum_{n\geq 0}R_n(G)$ is a graded commutative algebra over
${\Bbb Z}_2$ with unit. The multiplication in $R_*(G)$ is given by
$[V_1]\cdot[V_2]=[V_1\oplus V_2]$. We can identify $R_*(G)$ with the
graded polynomial algebra over ${\Bbb Z}_2$ generated by $\H(G,{\Bbb
Z}_2)$, where the addition in $\H(G,{\Bbb Z}_2)$ is given by the
tensor product of representations
$(\rho+\mu)(g)=\rho(g)\cdot\mu(g)$, and the multiplication is given
by the direct sum of representations. The homomorphisms
$\rho_i:(g_1,...g_n)\longmapsto g_i$ give a standard basis of
$\H(G,{\Bbb Z}_2)$. Then $R_*(G)$ is isomorphic to the graded
polynomial algebra ${\Bbb Z}_2[\rho_1,...,\rho_n]$. Obviously, each
$[\tau_pM]$ of $\mathcal{N}_M$ uniquely corresponds to a monomial of
degree $n$ in  ${\Bbb Z}_2[\rho_1,...,\rho_n]$ such that the $n$
factors of the monomial form a basis of $\H(G,{\Bbb Z}_2)$.

\vskip .2cm

There is a natural homomorphism $\delta_n:{\frak
M}_n\longrightarrow R_n(G)$ defined by
$$\delta_n([M,\phi])=\sum_{p\in M^G}[\tau_p M].$$
 The following result is essentially due to Stong \cite{s}.
\begin{thm}[Stong] \label{s}
$\delta_n$ is a monomorphism.
\end{thm}

Theorem~\ref{s} implies that for each $\beta$ in ${\frak M}_n$,
there must be a representative $(M,\phi)$ of $\beta$ such that
$\mathcal{N}_M$ is {\em prime} (i.e., either all elements of
$\mathcal{N}_{M}$ are distinct or $\mathcal{N}_{M}$ is empty), and
$\mathcal{N}_{M}$ is uniquely determined by $\beta$. Define
$$\mathcal{N}_{\beta}:=\mathcal{N}_{M}$$
and it is called the prime tangent representation set of $\beta$.
Then

\begin{cor} \label{ns}
Let $\beta_1,\beta_2\in {\frak M}_n$. Then
$$\beta_1=\beta_2\Longleftrightarrow
\mathcal{N}_{\beta_1}=\mathcal{N}_{\beta_2}.$$
\end{cor}

 \begin{rem} Since Hom$(G,{\Bbb Z}_2)$ is
isomorphic to $G$, each $[\tau_pM]$ of $\mathcal{N}_M$ actually
corresponds a unique element (denoted by $[\Delta_p]$) in the
quotient $\text{GL}(n,{\Bbb Z}_2)/\text{\bf S}_n$, where $\text{\bf
S}_n$ is a subgroup generated by all matrices of the form $E_{ij}$
(i.e., the identity matrix $E$ makes an exchange between the $i$-th
column and the $j$-th column), and it is isomorphic to the symmetric
group of rank $n$. Thus,  for any two $\sigma_1,\sigma_2$ in
$[\Delta_p]$, there exists a matrix $\theta$ in $\text{\bf S}_n$
such that $\sigma_1=\sigma_2\theta$. This also means that there is a
one-to-one correspondence between all bases of $({\Bbb Z}_2)^n$ and
$\text{GL}(n,{\Bbb Z}_2)/\text{\bf S}_n$. Here we call $[\Delta_p]$
the {\em tangent matrix} at $p$. With this understood,
 we often regard each element $[\tau_pM]$ of
$\mathcal{N}_M$ as being $[\Delta_p]$. Note that $|\text{GL}(n,{\Bbb
Z}_2)|=2^{{{n(n-1)}\over 2}}\prod_{i=1}^n(2^i-1)$, see \cite{ab}.
\end{rem}

\begin{prop} \label{bound}
Let $\beta$ be a nonzero element of ${\frak M}_n$. Then
$$n+1\leq \vert\mathcal{N}_\beta\vert\leq{{2^{{{n(n-1)}\over 2}}\prod_{i=1}^n(2^i-1)}\over {n!}}.$$
In particular, such upper and lower bounds are the best possible.
\end{prop}
\begin{proof}  The lower bound of $\vert
\mathcal{N}_\beta\vert$ is a special case of the Theorem 1.2 in
\cite{l}. Thus, it suffices to give the proof of the upper bound.
For this, it needs to merely show that there is a nonzero element
$\beta'\in {\frak M}_n$ such that $\vert
\mathcal{N}_{\beta'}\vert={{2^{{{n(n-1)}\over
2}}\prod_{i=1}^n(2^i-1)}\over {n!}}$. \vskip .2cm
 Consider the
standard $({\Bbb Z}_2)^n$-action $({\Bbb R}P^n, T_0)$ of $({\Bbb
Z}_2)^n$ on the real $n$-dimensional projective space ${\Bbb R}P^n$
defined by $n$ commuting involutions
$$t_i:([x_0, x_1,...,x_n])=[x_0, x_1,...,x_{i-1}, -x_i,x_{i+1},...,x_n], \ \ i=1,...,n$$
where $t_1,...,t_n$  generate  $({\Bbb Z}_2)^n$. This action fixes
$n+1$ isolated points
$$p_{i+1}=[\underbrace{0,...,0}_i, 1, 0,...,0]$$
where $i=0, 1,2,...,n$, and one easily sees that its tangent matrix
set is
\begin{eqnarray*} \mathcal{N}_0=\{ [\Delta_i]=\left[
\begin{pmatrix}
1 &  & & &  & & & & \\
   & \cdot &  & &  & & & & \\
    &    & \cdot &  & & & &  & \\
      &    &  & \cdot &  & &  & & \\
 1 &   & \cdots &  & 1 & & \cdots & & 1 \\
   &   &    &   &     & \cdot & &  & \\
  & &   &    &   &     & \cdot & & \\
 & & &   &    &   &    & \cdot  &  \\
& & & & &   &    &       &    1 \\
 \end{pmatrix}\right]
   \vert i=0,1,2,...,n\}
 \end{eqnarray*}
 where the row vector $(1,\cdots, 1,\cdots, 1)$ in $\Delta_i$ denotes
 $i$-th row, and one makes the convention that $\Delta_0=E$ when $i=0$;
 especially, each $[\Delta_i]$ corresponds to the isolated
 point $p_{i+1}$. Obviously, $\mathcal{N}_0$ is prime. By direct computations, one has that $\Delta_i\Delta_i=E$ and
 the result of the product $\Delta_i\Delta_j (i, j\not=0, j\not=i)$ just
 makes
 an exchange between $i$-th column and $j$-th column of
 $\Delta_j$. Thus, for $i, j\not=0$, one has
 \begin{equation} \label{r}
 [\Delta_i\Delta_j]
 =
 \begin{cases}
 [E] & \text{ if } i=j\\
 \text{$[\Delta_j]$}  & \text{ if } i\not=j.
 \end{cases}
 \end{equation}

 \vskip .2cm

Now, let $B_{n+1}$ denote the subset of $\text{GL}(n,{\Bbb Z}_2)$
defined as follows:
$$B_{n+1}=\{\sigma\in \text{GL}(n,{\Bbb Z}_2)\vert
\sigma\mathcal{N}_0=\mathcal{N}_0\}$$ where
$\sigma\mathcal{N}_0=\{[\sigma\Delta_0],[\sigma\Delta_1],...,[\sigma\Delta_n]\}$.
Obviously, $B_{n+1}$ is a subgroup of $\text{GL}(n,{\Bbb Z}_2)$,
and  each element of $B_{n+1}$ actually makes a permutation for
$[\Delta_0],[\Delta_1],...,[\Delta_n]$. One then knows from
(\ref{r}) that each $\Delta_i\in B_{n+1}$.

\vskip .2cm

{\bf Claim I.} $\vert B_{n+1}\vert=(n+1)!$.

\vskip .2cm

First, we prove that $B_{n+1}$ contains the symmetric group ${\bf
S}_n$. Actually, this is because for any $E_{ij}$ in ${\bf S}_n$
and  any $\Delta_l$,
 $$[E_{ij}\Delta_l]=\begin{cases}
 [\Delta_l] & \text{ if $i,j\not=l$ or $l=0$}\\
 [\Delta_i] & \text{ if $j=l\not=0$}\\
 [\Delta_j] & \text{ if $i=l\not=0$.}
 \end{cases}$$
Next, it is easy to see that for any $\sigma$ in $B_{n+1}$,
$\sigma$ can be expressed as a product by some matrices of ${\bf
S}_n$ and some of the $\Delta_i$'s.
 Obviously, when $i\not=0$, $\Delta_i\not\in {\bf S}_{n}$. Thus,
 $B_{n+1}$ is generated by those matrices of ${\bf S}_n$ and all
 $\Delta_i$, and Claim I then follows from this.

\vskip .2cm

{\bf Claim II.} For any $\sigma,\tau\in \text{GL}(n,{\Bbb Z}_2)$,
$$\sigma\mathcal{N}_0\cap\tau\mathcal{N}_0\not=\emptyset\Longleftrightarrow
\sigma\mathcal{N}_0=\tau\mathcal{N}_0.$$

\vskip .2cm

It is obvious that if $\sigma\mathcal{N}_0=\tau\mathcal{N}_0$ then
$\sigma\mathcal{N}_0\cap\tau\mathcal{N}_0\not=\emptyset$.
Conversely, if
$\sigma\mathcal{N}_0\cap\tau\mathcal{N}_0\not=\emptyset$, then
there are $[\Delta_i], [\Delta_j]\in \mathcal{N}_0$ such that
$[\sigma\Delta_j]=[\tau\Delta_i]$. By the definition of $B_{n+1}$,
one has that
$$\sigma\mathcal{N}_0=\tau\mathcal{N}_0\Longleftrightarrow
\sigma^{-1}\tau\in B_{n+1}.$$ Hence, it suffices to show that
$\sigma^{-1}\tau\in B_{n+1}$. From
$[\sigma\Delta_j]=[\tau\Delta_i]$,  there is an element $s\in
\text{\bf S}_n$ such that $\sigma\Delta_js=\tau\Delta_i$, so
$$\sigma^{-1}\tau=\Delta_j s\Delta_i.$$
Note that $\Delta_i^{-1}=\Delta_i$. Furthermore, one concludes
that $\sigma^{-1}\tau\in B_{n+1}$.

\vskip .2cm
 For  any  automorphism $\sigma:({\Bbb Z}_2)^n\longrightarrow
 ({\Bbb Z}_2)^n$ where $\sigma\in \text{GL}(n,{\Bbb Z}_2)$,
 one obtains new generators $\sigma (t_1),..., \sigma(t_n)$ of $({\Bbb Z}_2)^n$,
 and then one  obtains a new $({\Bbb Z}_2)^n$-action $({\Bbb R}P^n, \sigma T_0)$ from $({\Bbb
 R}P^n, T_0)$ by using generators $\sigma (t_1),..., \sigma(t_n)$ such
 that its tangent matrix set is $(\sigma^{-1})^\top\mathcal{N}_0$. By the
 above arguments with Corollary~\ref{ns} together, up to equivariant cobordism, there are
 ${{\vert\text{GL}(n,{\Bbb Z}_2)\vert}\over{\vert
 B_{n+1}\vert}}={{2^{{{n(n-1)}\over
2}}\prod_{i=1}^n(2^i-1)}\over {(n+1)!}}$ different $({\Bbb
Z}_2)^n$-actions $({\Bbb R}P^n, \sigma T_0)$, and especially,  the
union of their  tangent matrix
 sets just consists of  all elements of $\text{GL}(n,{\Bbb
 Z}_2)/{\bf S}_{n}$.
 Therefore, taking
 $$(T', {M'}^n)=\bigsqcup_{\{(\sigma^{-1})^\top\}\in\text{GL}(n,{\Bbb
 Z}_2)/B_{n+1}}({\Bbb R}P^n, \sigma T_0)$$
 then the tangent representation set of this action is prime, and the number of its all elements is
 $$(n+1)\times {{2^{{{n(n-1)}\over
2}}\prod_{i=1}^n(2^i-1)}\over {(n+1)!}}={{2^{{{n(n-1)}\over
2}}\prod_{i=1}^n(2^i-1)}\over {n!}}.$$ This completes the proof of
the upper bound.
\end{proof}

\section{Essential generators of ${\frak M}_n$}

\begin{defn} Let $\beta\not=0$ in ${\frak M}_n$. One says that $\beta$ is
{\em an essential generator} if
$\vert{\mathcal{N}}_{\beta+\gamma}\vert \geq
\vert{\mathcal{N}}_\beta\vert$ for any $\gamma\in {\frak M}_n$ with
$\vert{\mathcal{N}}_{\gamma}\vert < \vert{\mathcal{N}}_\beta\vert$.
\end{defn}

We know from Proposition~\ref{bound} that  up to equivariant
cobordism there are ${{2^{{{n(n-1)}\over
2}}\prod^n_{i=1}(2^i-1)}\over {(n+1)!}}$ different $({\Bbb
Z}_2)^n$-actions $({\Bbb R}P^n, \sigma T_0),
\sigma\in\text{GL}(n,{\Bbb Z}_2)$, and each $({\Bbb R}P^n, \sigma
T_0)$ fixes just $n+1$ isolated points with different
representations. Since the lower bound of
$\vert{\mathcal{N}}_\beta\vert$ for any nonzero element $\beta$ of
${\frak M}_n$ is $n+1$, we have that each $({\Bbb R}P^n, \sigma
T_0)$ is an essential generator.

\begin{lem} \label{l1}
Let $\beta\in {\frak M}_n$. If $\beta$ is an essential generator,
then
$$\vert{\mathcal{N}}_\beta\vert\leq
\begin{cases}
{{2^{{{n(n-1)}\over 2}}\prod^n_{i=1}(2^i-1)}\over {2n!}} & \text{
if $n$ is odd,}\\
  {{2^{{{n(n-1)}\over 2}}\prod^n_{i=1}(2^i-1)}\over {2(n-1)!(n+1)}}
  & \text{ if $n$ is even.}
  \end{cases}$$
\end{lem}
\begin{proof}
 If $\vert{\mathcal{N}}_\beta\vert=n+1$, then obviously the lemma holds.
  Now suppose that $\vert{\mathcal{N}}_\beta\vert>n+1$
so $\beta$ is not just one of the $({\Bbb R}P^n, \sigma T_0),
\sigma\in\text{GL}(n,{\Bbb Z}_2)$.  Then one claims that for each
$({\Bbb Z}_2)^n$-action $({\Bbb R}P^n, \sigma T_0)$,
${\mathcal{N}}_\beta$ cannot contain more than $[{{n+1}\over 2}]$
elements in ${\mathcal{N}}_{({\Bbb R}P^n, \sigma T_0)}$. Actually,
if not, then one has that
$\vert{\mathcal{N}}_\beta\vert>\vert{\mathcal{N}}_{\beta+[({\Bbb
R}P^n, \sigma T_0)]}\vert$, but this is impossible since $\beta$
is an essential generator.   The lemma then follows from this
claim.
\end{proof}
\begin{prop} \label{l2}
Let $\beta\in {\frak M}_n$. Then $\beta$ is a linear combination
of essential generators.
\end{prop}
\begin{proof}
The argument is trivial if $\beta=0$ or $\beta$ is an essential
generator. Suppose that $\beta$ is nonzero and is not an essential
generator. Then there exists some element $\gamma$ with
$\vert{\mathcal{N}}_{\gamma}\vert < \vert{\mathcal{N}}_\beta\vert$
in ${\frak M}_n$ such that $\beta=(\beta+\gamma)+\gamma$ with
$\vert{\mathcal{N}}_{\beta+\gamma}\vert <
\vert{\mathcal{N}}_\beta\vert$. If $\gamma$ or $\beta+\gamma$ is not
an essential generator, since ${\frak M}_n$ contains finite
elements, by continuing the above process, finally $\beta$ may be
expressed as  a linear combination of essential generators.
\end{proof}

\section{Small covers}

\vskip .2cm An $n$-dimensional convex polytope $P^n$ is said to be
{\em simple} if exactly $n$ faces of codimension one meet at each of
its vertices. Each point of a simple convex polytope $P^n$ has a
neighborhood which is affine isomorphic to an open subset of the
positive cone ${\Bbb R}_{\geq 0}^n$.  A smooth closed $n$-manifold
$M^n$ is said to be a {\em small cover} if it admits an effective
smooth $({\Bbb Z}_2)^n$-action and is locally isomorphic to the
standard action of $({\Bbb Z}_2)^n$ on ${\Bbb R}^n$ such that  the
orbit space of the action is a simple convex polytope $P^n$.

\vskip .2cm
 A small cover is a special  2-torus manifold. A
canonical example of small cover is the $n$-dimensional real
projective space ${\Bbb R}P^n$ with the standard $({\Bbb
Z}_2)^n$-action whose orbit space is the $n$-simplex $\Delta^n$.

\vskip .2cm Suppose that $\pi:M^n\longrightarrow P^n$ is a small
cover over a simple convex polytope $P^n$. Let
$\mathcal{F}(P^n)=\{F_1,...,F_\ell\}$ be the set of codimension-one
faces (facets) of $P^n$. Then there are $\ell$ connected
submanifolds $M_1,...,M_\ell$ determined by $\pi$ and $F_i$ (i.e.,
$M_i=\pi^{-1}(F_i)$), which are called {\em characteristic
submanifolds} here. Each submanifold $M_i$ is fixed pointwise by the
${\Bbb Z}_2$-subgroup $G_i$ of $({\Bbb Z}_2)^n$, so that each facet
$F_i$ corresponds to the ${\Bbb Z}_2$-subgroup $G_i$. Since there is
a canonical isomorphism from $({\Bbb Z}_2)^n$ to $\H({\Bbb Z}_2,
({\Bbb Z}_2)^n)$, such the ${\Bbb Z}_2$-subgroup $G_i$ corresponds
to an element $\upsilon_i$ in $\H({\Bbb Z}_2, ({\Bbb Z}_2)^n)$.
 For each face $F$ of
codimension $s$, since $P^n$ is simple, there are $s$ facets
$F_{i_1},...,F_{i_s}$  such that
$$F=F_{i_1}\cap\cdots \cap F_{i_s}.$$ Then, the corresponding
characteristic submanifolds $M_{i_1},...,M_{i_s}$ intersect
transversally in the $(n-s)$-dimensional submanifold
$\pi^{-1}(F)$, and the isotropy subgroup $G_F$ of $\pi^{-1}(F)$ is
a  subtorus of rank $s$ and is generated by $G_{i_1},...,G_{i_s}$
(or is determined by $\upsilon_{i_1},...,\upsilon_{i_s}$ in
$\H({\Bbb Z}_2, ({\Bbb Z}_2)^n)$).
 Thus,  this actually gives a characteristic
function (see \cite{dj})
$$\lambda:\mathcal{F}(P^n)\longrightarrow \H({\Bbb Z}_2, ({\Bbb Z}_2)^n)$$
defined by $\lambda(F_i)=\upsilon_i$ such that for any face
$F=F_{i_1}\cap\cdots \cap F_{i_s}$ of $P^n$,
$\lambda(F_{i_1}),...,\lambda(F_{i_s})$ are linearly independent in
$\H({\Bbb Z}_2, ({\Bbb Z}_2)^n)$. When $\dim F=0$ (i.e., $s=n$), $F$
is a vertex of $P^n$, which corresponds to a $({\Bbb Z}_2)^n$-fixed
point $p$ of $M$. In this case,
$\lambda(F_{i_1}),...,\lambda(F_{i_n})$ uniquely determines a dual
basis of $\H(({\Bbb Z}_2)^n, {\Bbb Z}_2)$, which just gives the
tangent  representation at $p$. Thus, the characteristic function
$\lambda$ completely determines the tangent representation set
$\mathcal{N}_M$ of fixed points of $M^n$.
 \vskip .2cm
By the work of Davis and Januszkiewicz \cite{dj}, there is a
reconstruction process of $M^n$ by using the product bundle $({\Bbb
Z}_2)^n\times P^n$ and $\lambda$. Note that each point $q\in
\partial P^n$ must lie in the relative interior of a unique face
$F(q)$ of $P^n$. Then, one may define an equivalence relation on
$({\Bbb Z}_2)^n\times P^n$ as follows:
$$(t_1, x)\sim (t_2, x)\Longleftrightarrow t_1^{-1}t_2\in
G_{F(q)}$$ where $x\in F(q)$, so that the quotient space
$$M(\lambda):=({\Bbb Z}_2)^n\times P^n/(t_1, x)\sim (t_2, x) $$
is equivariantly homeomorphic to $M^n$. Obviously, both $M^n$ and
$M(\lambda)$ have the same characteristic function, so they also are
cobordant equivariantly.

\vskip .2cm  By $\Lambda(P^n)$ we denote  the set of all
characteristic functions on $P^n$.   Then we have

\begin{prop}
 Let $\pi:M^n\longrightarrow P^n$ be a small cover over a simple
convex polytope $P^n$. Then all small covers over $P^n$ are given by
$\{M(\lambda)\vert \lambda\in \Lambda(P^n)\}$ from the viewpoint of
cobordism.
\end{prop}

  \noindent {\bf Remark.} Generally speaking, one cannot make sure
that there always exist characteristic functions (or colorings)
over a simple convex polytope $P^n$ when $n\geq 4$.
    For example, see [DJ, Nonexamples 1.22]. However,
    the Four Color Theorem makes sure that every 3-dimensional simple
 convex polytope  always admits characteristic functions.

 \vskip .2cm

The correspondence $\lambda\longmapsto\sigma\circ\lambda$ defines an
action of $\text{GL}(n,{\Bbb Z}_2)$ on $\Lambda(P^n)$, and it then
induces an action of $\text{GL}(n,{\Bbb Z}_2)$ on $\{M(\lambda)\vert
\lambda\in \Lambda(P^n)\}$, given by $M(\lambda)\longmapsto
M(\sigma\circ \lambda)$.  It is easy to check that such two actions
are free.

\vskip .2cm

 The following two kinds of small covers play a key important
role on indicating the structure of ${\frak M}_3$.

\begin{exam}[Small covers over a 3-complex $\Delta^3$] \label{e1}
{\em A 3-simplex $\Delta^3$ has four 2-faces, and a canonical
characteristic function $\lambda_0$ on it is defined by assigning to
$\rho_1^*,\rho_2^*,\rho_3^*,\rho_1^*+\rho_2^*+\rho_3^*$ the four
2-faces of $\Delta^3$, where $\{\rho_1^*,\rho_2^*,\rho_3^*\}$ is the
standard basis of $\H({\Bbb Z}_2, ({\Bbb Z}_2)^3)$, which
corresponds to $\rho_1,\rho_2,\rho_3$ of $\H(({\Bbb Z}_2)^3, {\Bbb
Z}_2)$. Thus, $\{\sigma\circ \lambda_0|\sigma\in\text{GL}(3,{\Bbb
Z}_2)\}$ gives all characteristic functions on $\Delta^3$. Since the
characteristic function of  the standard action $T_0$ of $({\Bbb
Z}_2)^3$ on ${\Bbb R}P^3$ is just $\lambda_0$, $\{M(\sigma\circ
\lambda_0)|\sigma\in\text{GL}(3,{\Bbb Z}_2)\}=\{({\Bbb R}P^3, \sigma
T_0)|\sigma\in\text{GL}(3,{\Bbb Z}_2)\}$. Proposition~\ref{bound}
has shown that, up to equivariant cobordism, there are 7 different
small covers in $\{M(\sigma\circ
\lambda_0)|\sigma\in\text{GL}(3,{\Bbb Z}_2)\}=\{({\Bbb R}P^3, \sigma
T_0)|\sigma\in\text{GL}(3,{\Bbb Z}_2)\}$, denoted by $({\Bbb R}P^3,
T_0), ({\Bbb R}P^3, T_1),..., ({\Bbb R}P^3, T_6)$, respectively. A
direct calculation gives the following table about the tangent
representation sets of seven different small covers.

\vskip .2cm \noindent
\begin{tabular}{|l|l|l|}
\multicolumn{2}{c}{Table I}\\[5pt]
 \hline \multicolumn{1}{|c|}{ Small cover
$M$} &
\multicolumn{1}{|c|}{tangent representation set $\mathcal{N}_M$} \\
\hline $({\Bbb R}P^3, T_0)$ & {\scriptsize $\rho_1\rho_2\rho_3,
\rho_1(\rho_1+\rho_2)(\rho_1+\rho_3),\rho_2(\rho_1+\rho_2)(\rho_2+\rho_3),\rho_3(\rho_1+\rho_3)(\rho_2+\rho_3)$}\\
\hline $({\Bbb R}P^3, T_1)$ & {\scriptsize
$\rho_1(\rho_1+\rho_2)(\rho_1+\rho_2+\rho_3),
\rho_1\rho_2(\rho_2+\rho_3),\rho_2\rho_3(\rho_1+\rho_2),\rho_3(\rho_2+\rho_3)(\rho_1+\rho_2+\rho_3)$}\\
\hline $({\Bbb R}P^3, T_2)$ & {\scriptsize
$\rho_1(\rho_1+\rho_3)(\rho_1+\rho_2+\rho_3),
\rho_1\rho_3(\rho_2+\rho_3),\rho_2\rho_3(\rho_1+\rho_3),\rho_2(\rho_2+\rho_3)(\rho_1+\rho_2+\rho_3)$}\\
\hline $({\Bbb R}P^3, T_3)$ & {\scriptsize
$\rho_2(\rho_1+\rho_2)(\rho_1+\rho_2+\rho_3),
\rho_1\rho_2(\rho_1+\rho_3),\rho_1\rho_3(\rho_1+\rho_2),\rho_3(\rho_1+\rho_3)(\rho_1+\rho_2+\rho_3)$}\\
\hline $({\Bbb R}P^3, T_4)$ & {\tiny    
$\rho_1(\rho_1+\rho_2)(\rho_2+\rho_3),
\rho_1\rho_2(\rho_1+\rho_2+\rho_3),\rho_2(\rho_1+\rho_2)(\rho_1+\rho_3),(\rho_1+\rho_3)
(\rho_2+\rho_3)(\rho_1+\rho_2+\rho_3)$}\\
\hline $({\Bbb R}P^3, T_5)$ & {\tiny    
$\rho_1(\rho_1+\rho_3)(\rho_2+\rho_3),
\rho_1\rho_3(\rho_1+\rho_2+\rho_3),\rho_3(\rho_1+\rho_2)(\rho_1+\rho_3),(\rho_1+\rho_2)
(\rho_2+\rho_3)(\rho_1+\rho_2+\rho_3)$}\\
\hline $({\Bbb R}P^3, T_6)$ & {\tiny    
$\rho_2(\rho_1+\rho_3)(\rho_2+\rho_3),
\rho_2\rho_3(\rho_1+\rho_2+\rho_3),(\rho_1+\rho_2)
(\rho_1+\rho_3)(\rho_1+\rho_2+\rho_3),\rho_3(\rho_1+\rho_2)(\rho_2+\rho_3)$}\\
\hline
\end{tabular} }
\end{exam}

\begin{exam}[Small covers over a prism  $P^3$]\label{e2}
{\em There exists only one simple convex 3-polytope with six
vertices (i.e., a prism $P^3$), see \cite{e}. Let $F_1, F_2, F_4$
denote three square facets, and $F_3, F_5$ two triangular facets in
$P^3$. From \cite{ccl}  we know that essentially there are five
different characteristic functions
$\lambda_1,\lambda_2,\lambda_3,\lambda_4,\lambda_5$ under the action
of $\text{GL}(3,{\Bbb Z}_2)$ on $\Lambda(P^3)$, which are
respectively defined by

$$\begin{tabular}{|l|l|l|l|l|l|}
 \hline  & $F_1$& $F_2$ & $F_3$ & $F_4$ & $F_5$
\\
\hline $\lambda_1$ & $\rho_1^*$ & $\rho_2^*$ & $\rho_3^*$ &
$\rho_1^*+\rho_2^*$ & $\rho_1^*+\rho_2^*+\rho_3^*$\\
\hline $\lambda_2$ & $\rho_1^*$ & $\rho_2^*$ & $\rho_3^*$ &
$\rho_1^*+\rho_2^*$ & $\rho_1^*+\rho_3^*$\\
\hline $\lambda_3$ & $\rho_1^*$ & $\rho_2^*$ & $\rho_3^*$ &
$\rho_1^*+\rho_2^*$ & $\rho_2^*+\rho_3^*$\\
\hline $\lambda_4$ & $\rho_1^*$ & $\rho_2^*$ & $\rho_3^*$ &
$\rho_1^*+\rho_2^*$ & $\rho_3^*$\\
\hline $\lambda_5$ & $\rho_1^*$ & $\rho_2^*$ & $\rho_3^*$ &
$\rho_1^*+\rho_2^*+\rho_3^*$ & $\rho_3^*$\\
\hline
\end{tabular}$$
It is easy to check that for any $\sigma\in \text{GL}(3,{\Bbb
Z}_2)$,  every one of $M(\sigma\circ \lambda_4)$ and $M(\sigma\circ
\lambda_5)$ always bounds equivariantly. A direct calculation
shows that for $\sigma_1=\begin{pmatrix} 1 & & \\
  1  & 1 & \\
    & & 1
    \end{pmatrix}$,
    $\mathcal{N}_{M(\sigma_1\circ \lambda_1)}=\mathcal{N}_{M(\lambda_2)}$,
    and for $\sigma_2=
    \begin{pmatrix} 1 &  1& \\
   & 1 & \\
    & & 1
    \end{pmatrix}$,
    $\mathcal{N}_{M(\sigma_2\circ \lambda_1)}=\mathcal{N}_{M(\lambda_3)}$.
    Since $\mathcal{N}_{M(\lambda_1)}$ is prime, by  Corollary~\ref{ns}, all nonzero equivariant cobordism classes in
    $\{M(\sigma\circ\lambda_1)|\sigma\in \text{GL}(3,{\Bbb
Z}_2)\}$ give those in all small covers over $P^3$. By further
computations, one obtains that there are only four matrices
$$\tau_1=\begin{pmatrix} 1 & & \\
    & 1 & \\
    & & 1
    \end{pmatrix}, \tau_2=\begin{pmatrix} 1 & & \\
    1& 1 & \\
    & & 1
    \end{pmatrix},\tau_3=\begin{pmatrix} 1 & & \\
    & 1 & 1\\
    & & 1
    \end{pmatrix},\tau_4=\begin{pmatrix} 1 & & \\
    1& 1 &1 \\
    & & 1
    \end{pmatrix}$$
such that
$\tau_i\mathcal{N}_{M(\lambda_1)}=\mathcal{N}_{M(\lambda_1)},
i=1,2,3,4$, and these four matrices form a subgroup of
$\text{GL}(3,{\Bbb Z}_2)$. Thus, up to equivariant cobordism, there
are ${{|\text{GL}(3,{\Bbb Z}_2)\vert}\over 4}=42$ different
nonbouding small covers over $P^3$. We can even construct such small
covers as follows. Consider the $({\Bbb Z}_2)^3$-action $\Phi_0$ on
$S^1\times {\Bbb R}P^2=S^1\times {\Bbb R}P({\Bbb C}\oplus{\Bbb R})$
defined by the following three commutative involutions
$$t_1: (z, [v, w])\longmapsto (\bar{z}, [\bar{z}v, w])$$
$$t_2: (z, [v, w])\longmapsto (z, [z\bar{v}, w])$$
$$t_3: (z, [v, w])\longmapsto (z, [-z\bar{v}, w]).$$
This action fixes six isolated points $(\pm 1, [0, 1]), (\pm 1, [1,
0]), (\pm 1, [\sqrt{-1}, 0])$, and its orbit space is just a prime
$P^3$. A direct calculation shows that $\mathcal{N}_{(S^1\times
{\Bbb R}P^2, \Phi_0)}$ consists of six distinct monomials
$\rho_1\rho_2\rho_3,
\rho_1\rho_2(\rho_2+\rho_3),\rho_1\rho_3(\rho_2+\rho_3),\rho_1(\rho_1+\rho_2)(\rho_1+\rho_3),
\rho_1(\rho_1+\rho_2)(\rho_2+\rho_3),\rho_1(\rho_1+\rho_3)(\rho_2+\rho_3)$
of ${\Bbb Z}_2[\rho_1,\rho_2,\rho_3]$, so $(S^1\times {\Bbb R}P^2,
\Phi_0)$ is nonbounding. Further, up to equivariant cobordism, 42
different nonbouding small covers over $P^3$ can be given by
applying automorphisms of $({\Bbb Z}_2)^3$ to $(S^1\times {\Bbb
R}P^2, \Phi_0)$, and they are denoted by $(S^1\times {\Bbb R}P^2,
\Phi_0)$, $(S^1\times {\Bbb R}P^2, \Phi_1)$, ..., $(S^1\times {\Bbb
R}P^2, \Phi_{41})$, respectively. }
\end{exam}
\vskip .2cm

\section{Graphs of actions}

 Given a nonzero element $\beta$ in
${\frak M}_n$. Let $(M^n,\phi)$ be a representative of $\beta$ such
that ${\mathcal{N}}_M$ is prime. Choose a nontrivial irreducible
representation $\rho$ in $\H(({\Bbb Z}_2)^n,{\Bbb Z}_2)$, let $C$ be
a component of the fixed point set of $\ker\rho(\cong ({\Bbb
Z}_2)^{n-1})$ acting on $M$ such that $\dim C>0$, and  the action of
$({\Bbb Z}_2)^n/\ker\rho$ on $C$ has a nonempty fixed point set.
Then the dimension of $C$ must be 1 since the action is effective,
and thus $C$ is equivariantly diffeomorphic to the circle $S^1$ with
a reflection fixing just two fixed points. Then one has an edge
joining these two fixed points, which is labeled by $\rho$.
Furthermore, one can obtain a graph $\Gamma_M$, which is the union
of all those edges chosen for each $\rho$ and $C$. Clearly, the set
of vertices of $\Gamma_M$ is just the fixed point set of $({\Bbb
Z}_2)^n$ acting on $M$. Since the tangent representation at a fixed
point $p$ has $n$ irreducible summands, the number of edges  in
$\Gamma_M$ meeting at $p$ is exactly $n$, so $\Gamma_M$ is a regular
graph of valence $n$. It should be pointed out that, generally,
$\Gamma_M$ is not determined by $\beta$ uniquely, and it depends
upon the choice of representatives of $\beta$.

\vskip .2cm

Let $E_{\Gamma_M}$ denote the set of all edges in $\Gamma_M$, and
let $V_{\Gamma_M}$ denote the set of all vertices in $\Gamma_M$.
Given a vertex $p$ in $V_{\Gamma_M}$, let $E_p$ denote the set of
$n$ edges joining to $p$. Then there is a natural map $\alpha:
E_{\Gamma_M}\longrightarrow \H(({\Bbb Z}_2)^n, {\Bbb Z}_2)$ (called
an axial function or a $({\Bbb Z}_2)^n$-coloring, cf~\cite{gz1},
\cite{gz2}, \cite{bl}). One knows from \cite{l} that $\alpha$
satisfies the following properties:

\vskip .2cm

1) for each vertex $p$ in $V_{\Gamma_\beta}$, $\alpha(E_p)$ spans
$\H(({\Bbb Z}_2)^n, {\Bbb Z}_2)$;

2) for each edge $e$ in $E_{\Gamma_\beta}$,
$$\prod_{x\in E_p-E_e}\alpha(x)\equiv \prod_{y\in E_q-E_e}\alpha(y)\mod \alpha(e)$$
where $p, q$ are two endpoints of $e$, and $E_e$ denotes the set of
all edges joining two endpoints of $e$. The pair $(\Gamma_M,
\alpha)$ is called the {\em moment graph} of $(M^n,\phi)$. Since
${\mathcal{N}}_M$ is prime, one has from \cite{l} that for each edge
$e$ in $\Gamma_M$, $\vert E_e\vert=1$.

\vskip .2cm

{\em Note.} If $M$ is a small cover over a simple convex polytope
$P^n$, then  $\Gamma_M$ is just the 1-skeleton of $P^n$. In this
case, it is easy to see that the map $\alpha:
E_{\Gamma_M}\longrightarrow \H(({\Bbb Z}_2)^n, {\Bbb Z}_2)$ is dual
to the characteristic function
$\lambda:\mathcal{F}(P^n)\longrightarrow \H({\Bbb Z}_2,({\Bbb
Z}_2)^n)$. In other words, both $\alpha$ and $\lambda$ are
determined to each other.

\vskip .2cm

By \cite{bl} we know that $(\Gamma_M,\alpha)$ is a ``good'' $({\Bbb
Z}_2)^n$-coloring, so that each $k$-nest $\Delta^k$ of
$(\Gamma_M,\alpha)$ is a connected regular $k$-valent subgraph of
$\Gamma_M$ with $\dim\text{Span}\alpha(\Delta^k)=k$, where
$\text{Span}\alpha(\Delta^k)$ denotes the linear space spanned by
all colors of edges in $\Delta^k$.   By $\mathcal{K}_{(\Gamma_M,
\alpha)}$ one denotes the set of all nests of $(\Gamma_M, \alpha)$.
Since each $k$-nest $(k>0)$ determines a $k$-dimensional subspace of
$\H(({\Bbb Z}_2)^n, {\Bbb Z}_2)$, it corresponds to an
$(n-k)$-dimensional subspace in the dual space $\H({\Bbb Z}_2,
({\Bbb Z}_2)^n)$. This actually gives a dual map $\eta$ from
$\mathcal{K}_{(\Gamma_M, \alpha)}$ to the set of all subspaces of
$\H({\Bbb Z}_2, ({\Bbb Z}_2)^n)$, which is just the characteristic
function when $M$ is a small cover. Obviously, $\eta$ maps each
$(n-1)$-dimensional nest  of $\mathcal{K}_{(\Gamma_M, \alpha)}$ to a
nonzero element in $\H({\Bbb Z}_2,({\Bbb Z}_2)^n)$. Since each
vertex $p$ is the intersection of $n$ $(n-1)$-nests of
$\mathcal{K}_{(\Gamma_M, \alpha)}$, it corresponds to a basis of
$\H({\Bbb Z}_2,({\Bbb Z}_2)^n)$, which is just the dual basis of the
basis $\alpha(E_p)$ in $\H(({\Bbb Z}_2)^n,{\Bbb Z}_2)$.

\vskip .2cm

From \cite{bl}  one knows that if $\dim M\leq 3$, then $(\Gamma_M,
\alpha)$ always admits a skeletal expansion (note that if $\dim
M>4$, under what condition $(\Gamma_M, \alpha)$ admits a skeletal
expansion  is still open). This will directly lead us to use this
result to study  the group structure of ${\frak M}_3$.

\begin{prop} [\cite{bl}]\label{gr}
If $\dim M=3$, then  $(\Gamma_M,\alpha)$  admits a 2-skeletal
expansion $(N,K)$ such that $N$ is a closed surface.
\end{prop}

 \section{Determination of ${\frak M}_3$}

The main task of this section is  devoted to determining the
structure of ${\frak M}_3$.

\begin{lem}\label{l4}
Let $\beta\in {\frak M}_3$. Then $\vert{\mathcal{N}}_\beta\vert$
is even.
\end{lem}
\begin{proof} The Euler characteristic of any 3-dimensional closed
manifold is always zero, and the lemma then follows from the
classical Smith Theorem.
\end{proof}

By Proposition~\ref{l2}, the key point of determining the structure
of ${\frak M}_3$ is to find out all essential generators in ${\frak
M}_3$. The following proposition characterizes the essential
generators of ${\frak M}_3$.

\begin{prop} \label{p}
A nonzero element $\beta\in{\frak M}_3$ is an essential generator
if and only if $\vert{\mathcal{N}}_\beta\vert\leq 6$. Further, all
essential generators of ${\frak M}_3$ are given by two kinds of
small covers $({\Bbb R}P^3, \sigma T_0)$ and $(S^1\times {\Bbb
R}P^2, \sigma \Phi_0)$.
\end{prop}

\begin{lem}\label{l5}
Let $\beta\in{\frak M}_3$ be  nonzero. If $\vert{\mathcal{N}}_\beta\vert\leq
6$, then $\beta$ is an essential generator.
\end{lem}
\begin{proof} If $\vert{\mathcal{N}}_\beta\vert=4$, then $\beta$
is one of $[({\Bbb R}P^3, \sigma T_0)]$'s so $\beta$ is an essential
generator. Thus, it suffices to consider the case
$\vert{\mathcal{N}}_\beta\vert=6$ by Lemma~\ref{l4}.  From
Example~\ref{e1}, we see that all $\mathcal{N}_{({\Bbb R}P^3, T_i)},
i=0,1,...,6$ are disjoint to each other. We first claim that any
intersection $\mathcal{N}_\beta\cap \mathcal{N}_{({\Bbb R}P^3,
T_i)}$ cannot contain four elements. If not, then there exists  some
$i^{'}$ such that $|\mathcal{N}_{\beta+[({\Bbb R}P^3,
T_{i^{'}})]}|=2$. By \cite{ks}, $\beta+[({\Bbb R}P^3, T_{i^{'}})]$
must be zero in ${\frak M}_3$, but this is impossible. Next, we
shall prove that any intersection $\mathcal{N}_\beta\cap
\mathcal{N}_{({\Bbb R}P^3, T_i)}$ cannot contain three elements. If
not, then there exists  some $i^{''}$ such that
$|\mathcal{N}_{\beta+[({\Bbb R}P^3, T_{i^{''}})]}|=4$, so that
$\beta+[({\Bbb R}P^3, T_{i^{''}})]$ must be the equivariant
cobordism class of another $({\Bbb R}P^3, T_j)$ with $j\not=i^{''}$.
Further, $\beta$ is the sum $[({\Bbb R}P^3, T_{i^{''}})]+[({\Bbb
R}P^3, T_j)]$, so $|\mathcal{N}_\beta|$ is 8 rather than 6. This is
a contradiction. Combining the above argument, one has that
$|\mathcal{N}_\beta\cap \mathcal{N}_{({\Bbb R}P^3, T_i)}|$ is less
than 3. Then the lemma follows from this.
\end{proof}

The following lemma indicates the connection between
$\mathcal{A}=\{({\Bbb R}P^3, T_i)| i=0,1,...,6\}$ and
$\mathcal{B}=\{(S^1\times {\Bbb R}P^2, \Phi_j)| j=0,1,...,41\}$.

\begin{lem}\label{l6}
Each $({\Bbb R}P^3, T_i)$ of $\mathcal{A}$  corresponds to six
small covers $(S^1\times {\Bbb R}P^2, \Phi_{i_1})$,...,
$(S^1\times {\Bbb R}P^2, \Phi_{i_6})$ of $\mathcal{B}$ for which
such six small covers are not cobordant to each other, and
$|\mathcal{N}_{({\Bbb R}P^3, T_i)}\cap \mathcal{N}_{(S^1\times
{\Bbb R}P^2, \Phi_{i_u})}|=2, u=1,...,6$.
\end{lem}

\begin{proof}
Since all $\mathcal{N}_{({\Bbb R}P^3, T_i)}, i=0,1,...,6$, are
distinct and since all $({\Bbb R}P^3, T_i), i=0,1,...,6$, can be
translated to each other up to cobordism by applying automorphisms
of $({\Bbb Z}_2)^3$, it suffices to consider the case of $({\Bbb
R}P^3, T_0)$. We see from the table I of Example~\ref{e1} that
$$\mathcal{N}_{({\Bbb R}P^3, T_0)}=\{\rho_1\rho_2\rho_3,
\rho_1(\rho_1+\rho_2)(\rho_1+\rho_3),\rho_2(\rho_1+\rho_2)(\rho_2+\rho_3),\rho_3(\rho_1+\rho_3)(\rho_2+\rho_3)\}.$$
Obviously, any two monomials  of $\mathcal{N}_{({\Bbb R}P^3, T_0)}$
give five elements of $\H(({\Bbb Z}_2)^3,{\Bbb Z}_2)$, and there are
exactly six such pairs in $\mathcal{N}_{({\Bbb R}P^3, T_0)}$.
Consider two monomials $\rho_1\rho_2\rho_3,
\rho_1(\rho_1+\rho_2)(\rho_1+\rho_3)$ of $\mathcal{N}_{({\Bbb R}P^3,
T_0)}$, we get five elements $\rho_1, \rho_2, \rho_3, \rho_1+\rho_2,
\rho_1+\rho_3$ of $\H(({\Bbb Z}_2)^3,{\Bbb Z}_2)$. Using these five
elements, we can define an axial function $\alpha$ on the 1-skeleton
of a prism $P^3$ as shown in Figure~\ref{a1}.
\begin{figure}[h]
    \input{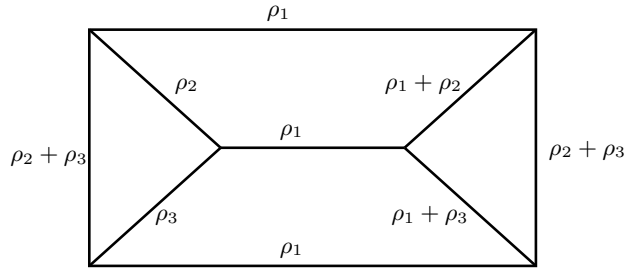}\centering
    \caption[]{ An axial function
$\alpha$ on the 1-skeleton of a prism $P^3$}\label{a1}
\end{figure}
Since $\alpha$ uniquely determines a characteristic function on
$P^3$, we obtain a small cover $(S^1\times {\Bbb R}P^2, \Phi_{0_1})$
with six fixed points over $P^3$ such that its tangent
representation set $\mathcal{N}_{(S^1\times {\Bbb R}P^2,
\Phi_{0_1})}$ consists of six monomials $\rho_1\rho_2\rho_3,
\rho_1(\rho_1+\rho_2)(\rho_1+\rho_3), \rho_1\rho_2(\rho_2+\rho_3),
\rho_1\rho_3(\rho_2+\rho_3),\rho_1(\rho_1+\rho_2)(\rho_2+\rho_3),
\rho_1(\rho_1+\rho_3)(\rho_2+\rho_3)$. Similarly, for other five
pairs in $\mathcal{N}_{({\Bbb R}P^3, T_0)}$, we can obtain five
small covers $(S^1\times {\Bbb R}P^2, \Phi_{0_u}), u=2,...,6$ with
their tangent representation sets as follows {\tiny
$$\begin{tabular}{|l|l|}
\hline $u$ & \ \ \ \ \ \ \ \ \ \ \ \ \ \ \ \ \ \ \ \ \ \ \ \ \ \ \ \
\ \ \ \ \ \ \ \ \ \ \ \ \ \ \ \ \ \ \ \ \ \ \ \
$\mathcal{N}_{(S^1\times
{\Bbb R}P^2, \Phi_{0_u})}$\\
\hline 2 & $\{\rho_1\rho_2\rho_3,\rho_1\rho_2(\rho_1+\rho_3),
\rho_2\rho_3(\rho_1+\rho_3),\rho_2(\rho_1+\rho_2)(\rho_2+\rho_3),\rho_2(\rho_1+\rho_2)(\rho_1+\rho_3),
\rho_2(\rho_1+\rho_3)(\rho_2+\rho_3)\}$\\
\hline 3 & $\{\rho_1\rho_2\rho_3,\rho_1\rho_3(\rho_1+\rho_2),
\rho_2\rho_3(\rho_1+\rho_2),\rho_3(\rho_1+\rho_3)(\rho_2+\rho_3),\rho_3(\rho_1+\rho_2)(\rho_1+\rho_3),
\rho_3(\rho_1+\rho_2)(\rho_2+\rho_3)\}$\\
\hline 4 &
$\{\rho_1(\rho_1+\rho_2)(\rho_1+\rho_3),\rho_1\rho_3(\rho_1+\rho_2),
\rho_3(\rho_1+\rho_2)(\rho_1+\rho_3),\rho_2(\rho_1+\rho_2)(\rho_2+\rho_3),\rho_2\rho_3(\rho_1+\rho_2),
\rho_3(\rho_1+\rho_2)(\rho_2+\rho_3)\}$\\
\hline 5 &
$\{\rho_1(\rho_1+\rho_2)(\rho_1+\rho_3),\rho_1\rho_2(\rho_1+\rho_3),
\rho_2(\rho_1+\rho_2)(\rho_1+\rho_3),\rho_3(\rho_1+\rho_3)(\rho_2+\rho_3),\rho_2\rho_3(\rho_1+\rho_3),
\rho_2(\rho_1+\rho_3)(\rho_2+\rho_3)\}$\\
\hline 6 &
$\{\rho_2(\rho_1+\rho_2)(\rho_2+\rho_3),\rho_1\rho_2(\rho_2+\rho_3),
\rho_1(\rho_1+\rho_2)(\rho_2+\rho_3),\rho_3(\rho_1+\rho_3)(\rho_2+\rho_3),\rho_1\rho_3(\rho_2+\rho_3),
\rho_1(\rho_1+\rho_3)(\rho_2+\rho_3)\}$\\
\hline \end{tabular}$$ } \noindent  Then the lemma follows from the
above argument and Corollary~\ref{ns}.
\end{proof}

\vskip .2cm

\begin{rem}\label{re}
Lemma~\ref{l6} also gives the method of constructing 42 different
small covers (up to equivariant cobordism) with 6 fixed points. In
particular, we easily see the following property that for each
$(S^1\times {\Bbb R}P^2, \Phi_j)$, two of $\mathcal{N}_{(S^1\times
{\Bbb R}P^2, \Phi_j)}$ are in some $\mathcal{N}_{({\Bbb R}P^3,
T_i)}$, and others are just distributed in four different
$\mathcal{N}_{({\Bbb R}P^3, T_{i_1})}, \mathcal{N}_{({\Bbb R}P^3,
T_{i_2})}, \mathcal{N}_{({\Bbb R}P^3, T_{i_3})},
\mathcal{N}_{({\Bbb R}P^3, T_{i_4})}$ with $i_v\not=i, v=1,2,3,4$.
In addition, we also see from the argument of Lemma~\ref{l6} that
$$\delta_3([(S^1\times
{\Bbb R}P^2, \Phi_{0_1})]+[(S^1\times {\Bbb R}P^2,
\Phi_{0_6})]+[({\Bbb R}P^3, T_0)])=0$$
$$\delta_3([(S^1\times
{\Bbb R}P^2, \Phi_{0_2})]+[(S^1\times {\Bbb R}P^2,
\Phi_{0_5})]+[({\Bbb R}P^3, T_0)])=0$$
$$\delta_3([(S^1\times
{\Bbb R}P^2, \Phi_{0_3})]+[(S^1\times {\Bbb R}P^2,
\Phi_{0_4})]+[({\Bbb R}P^3, T_0)])=0$$ where $\delta_3$ is the
monomorphism of Theorem~\ref{s}. This means that actually we need
only to consider the half of  42 different small covers $(S^1\times
{\Bbb R}P^2, \Phi_j), j=0,1,...,41$, such that up to equivariant
cobordism the union of any two of them is not one of $({\Bbb R}P^3,
T_i), i=0,1,...,6$. With no loss we may assume that such 21
different small covers are just $(S^1\times {\Bbb R}P^2, \Phi_j),
j=0,1,...,20$, with their tangent representation sets stated in
Table II.
\end{rem}

 \noindent
\begin{tabular}{|l|l|l|}
\multicolumn{2}{c}{Table II}\\[5pt]
 \hline \multicolumn{1}{|c|}{ Small cover
$M$} &
\multicolumn{1}{|c|}{tangent representation set $\mathcal{N}_M$} \\
\hline {\scriptsize $(S^1\times {\Bbb R}P^2, \Phi_0)$} &
{\scriptsize $\rho_1\rho_2\rho_3,\rho_1\rho_2(\rho_2+\rho_3),
\rho_1\rho_3(\rho_2+\rho_3),\rho_1(\rho_1+\rho_2)(\rho_1+\rho_3),
\rho_1(\rho_1+\rho_2)(\rho_2+\rho_3),$}\\
&{\scriptsize $
\rho_1(\rho_1+\rho_3)(\rho_2+\rho_3)$}\\
\hline {\scriptsize $(S^1\times {\Bbb R}P^2, \Phi_1)$} &
{\scriptsize $\rho_1\rho_2\rho_3, \rho_1\rho_2(\rho_1+\rho_3),
\rho_2\rho_3(\rho_1+\rho_3), \rho_2(\rho_1+\rho_2)(\rho_1+\rho_3),
\rho_2(\rho_1+\rho_2)(\rho_2+\rho_3), $}\\ & {\scriptsize
$\rho_2(\rho_1+\rho_3)(\rho_2+\rho_3)$}\\
\hline {\scriptsize $(S^1\times {\Bbb R}P^2, \Phi_2)$} &
{\scriptsize $\rho_1\rho_2\rho_3, \rho_1\rho_3(\rho_1+\rho_2),
\rho_2\rho_3(\rho_1+\rho_2), \rho_3(\rho_1+\rho_2)(\rho_1+\rho_3),
\rho_3(\rho_1+\rho_2)(\rho_2+\rho_3),$}\\ & {\scriptsize
 $ \rho_3(\rho_1+\rho_3)(\rho_2+\rho_3)$}\\
\hline {\scriptsize $(S^1\times {\Bbb R}P^2, \Phi_3)$} &
{\scriptsize $\rho_1\rho_2\rho_3, \rho_1\rho_2(\rho_1+\rho_3),
\rho_2\rho_3(\rho_1+\rho_2), \rho_2\rho_3(\rho_1+\rho_3),
\rho_2\rho_3(\rho_1+\rho_2+\rho_3),$}\\ & {\scriptsize
 $ \rho_2(\rho_1+\rho_2)(\rho_1+\rho_2+\rho_3)$}\\
 \hline {\scriptsize $(S^1\times {\Bbb R}P^2, \Phi_4)$} &
{\scriptsize $\rho_1\rho_2(\rho_1+\rho_2+\rho_3),
\rho_1\rho_3(\rho_1+\rho_2+\rho_3),
\rho_1(\rho_1+\rho_2)(\rho_1+\rho_2+\rho_3),
 $}\\ & {\scriptsize
 $\rho_1(\rho_1+\rho_3)(\rho_1+\rho_2+\rho_3),
 \rho_2(\rho_1+\rho_2)(\rho_1+\rho_2+\rho_3), \rho_3(\rho_1+\rho_3)(\rho_1+\rho_2+\rho_3)$}\\
 \hline {\scriptsize $(S^1\times {\Bbb R}P^2, \Phi_5)$} &
{\scriptsize $\rho_1\rho_3(\rho_1+\rho_2),
\rho_1(\rho_1+\rho_2)(\rho_1+\rho_3),
\rho_2(\rho_1+\rho_2)(\rho_1+\rho_3),
 $}\\ & {\scriptsize
 $\rho_3(\rho_1+\rho_2)(\rho_1+\rho_3),
 \rho_2(\rho_1+\rho_2)(\rho_1+\rho_2+\rho_3), (\rho_1+\rho_2)(\rho_1+\rho_3)(\rho_1+\rho_2+\rho_3)$}\\
 \hline {\scriptsize $(S^1\times {\Bbb R}P^2, \Phi_6)$} &
{\scriptsize $\rho_1\rho_2\rho_3, \rho_1\rho_2(\rho_2+\rho_3),
\rho_1\rho_3(\rho_1+\rho_2),
 $}\\ & {\scriptsize
 $\rho_1\rho_3(\rho_2+\rho_3),
 \rho_1\rho_3)(\rho_1+\rho_2+\rho_3), \rho_1(\rho_1+\rho_2)(\rho_1+\rho_2+\rho_3)$}\\
 \hline {\scriptsize $(S^1\times {\Bbb R}P^2, \Phi_7)$} &
{\scriptsize $\rho_1\rho_2(\rho_1+\rho_2+\rho_3),
\rho_2\rho_3(\rho_1+\rho_2+\rho_3),
\rho_1(\rho_1+\rho_2)(\rho_1+\rho_2+\rho_3),
 $}\\ & {\scriptsize
 $\rho_2(\rho_1+\rho_2)(\rho_1+\rho_2+\rho_3),
 \rho_2(\rho_2+\rho_3)(\rho_1+\rho_2+\rho_3), \rho_3(\rho_2+\rho_3)(\rho_1+\rho_2+\rho_3)$}\\
 \hline {\scriptsize $(S^1\times {\Bbb R}P^2, \Phi_8)$} &
{\scriptsize $\rho_2\rho_3(\rho_1+\rho_2),
\rho_2(\rho_1+\rho_2)(\rho_2+\rho_3),
\rho_1(\rho_1+\rho_2)(\rho_2+\rho_3),
 $}\\ & {\scriptsize
 $\rho_3(\rho_1+\rho_2)(\rho_2+\rho_3),
 \rho_1(\rho_1+\rho_2)(\rho_1+\rho_2+\rho_3), (\rho_1+\rho_2)(\rho_2+\rho_3)(\rho_1+\rho_2+\rho_3)$}\\
 \hline {\scriptsize $(S^1\times {\Bbb R}P^2, \Phi_9)$} &
{\scriptsize $\rho_1\rho_2\rho_3, \rho_1\rho_2(\rho_1+\rho_3),
\rho_1\rho_2(\rho_2+\rho_3),
 $}\\ & {\scriptsize
 $\rho_1\rho_2(\rho_1+\rho_2+\rho_3),
 \rho_1\rho_3(\rho_2+\rho_3), \rho_1(\rho_1+\rho_3)(\rho_1+\rho_2+\rho_3)$}\\
 \hline
\end{tabular}

\noindent
\begin{tabular}{|l|l|l|}
 \hline {\scriptsize $(S^1\times {\Bbb R}P^2, \Phi_{10})$} &
{\scriptsize $\rho_1\rho_3(\rho_1+\rho_2+\rho_3),
\rho_2\rho_3(\rho_1+\rho_2+\rho_3),
\rho_1(\rho_1+\rho_3)(\rho_1+\rho_2+\rho_3),
 $}\\ & {\scriptsize
 $\rho_2(\rho_2+\rho_3)(\rho_1+\rho_2+\rho_3),
 \rho_3(\rho_1+\rho_3)(\rho_1+\rho_2+\rho_3), \rho_3(\rho_2+\rho_3)(\rho_1+\rho_2+\rho_3)$}\\
 \hline {\scriptsize $(S^1\times {\Bbb R}P^2, \Phi_{11})$} &
{\scriptsize $\rho_2\rho_3(\rho_1+\rho_3),
\rho_1(\rho_1+\rho_3)(\rho_2+\rho_3),
\rho_1(\rho_1+\rho_3)(\rho_1+\rho_2+\rho_3),
 $}\\ & {\scriptsize
 $\rho_2(\rho_1+\rho_3)(\rho_2+\rho_3),
 \rho_3(\rho_1+\rho_3)(\rho_2+\rho_3), (\rho_1+\rho_3)(\rho_2+\rho_3)(\rho_1+\rho_2+\rho_3)$}\\
 \hline {\scriptsize $(S^1\times {\Bbb R}P^2, \Phi_{12})$} &
{\scriptsize $\rho_1\rho_2(\rho_2+\rho_3),
\rho_1(\rho_1+\rho_2)(\rho_2+\rho_3),
\rho_2(\rho_1+\rho_2)(\rho_2+\rho_3),
 $}\\ & {\scriptsize
 $\rho_2(\rho_1+\rho_3)(\rho_2+\rho_3),
 \rho_2(\rho_2+\rho_3)(\rho_1+\rho_2+\rho_3), (\rho_1+\rho_3)(\rho_2+\rho_3)(\rho_1+\rho_2+\rho_3)$}\\
 \hline {\scriptsize $(S^1\times {\Bbb R}P^2, \Phi_{13})$} &
{\scriptsize $\rho_1\rho_2(\rho_1+\rho_3),
\rho_1\rho_2(\rho_1+\rho_2+\rho_3),
\rho_1(\rho_1+\rho_2)(\rho_1+\rho_3),
 $}\\ & {\scriptsize
 $\rho_1(\rho_1+\rho_2)(\rho_2+\rho_3),
 \rho_1(\rho_1+\rho_3)(\rho_2+\rho_3), \rho_1(\rho_1+\rho_3)(\rho_1+\rho_2+\rho_3)$}\\
 \hline {\scriptsize $(S^1\times {\Bbb R}P^2, \Phi_{14})$} &
{\scriptsize $\rho_1(\rho_1+\rho_2)(\rho_2+\rho_3),
\rho_1(\rho_1+\rho_2)(\rho_1+\rho_2+\rho_3),
\rho_2(\rho_1+\rho_2)(\rho_1+\rho_3),
 $}\\ & {\scriptsize
 $\rho_2(\rho_1+\rho_2)(\rho_1+\rho_2+\rho_3),
 (\rho_1+\rho_2)(\rho_1+\rho_3)(\rho_1+\rho_2+\rho_3), (\rho_1+\rho_2)(\rho_2+\rho_3)(\rho_1+\rho_2+\rho_3)$}\\
 \hline {\scriptsize $(S^1\times {\Bbb R}P^2, \Phi_{15})$} &
{\scriptsize $\rho_1\rho_3(\rho_2+\rho_3),
\rho_1(\rho_1+\rho_3)(\rho_2+\rho_3),
\rho_3(\rho_1+\rho_2)(\rho_2+\rho_3),
 $}\\ & {\scriptsize
 $\rho_3(\rho_1+\rho_3)(\rho_2+\rho_3),
 \rho_3(\rho_2+\rho_3)(\rho_1+\rho_2+\rho_3), (\rho_1+\rho_2)(\rho_2+\rho_3)(\rho_1+\rho_2+\rho_3)$}\\
 \hline {\scriptsize $(S^1\times {\Bbb R}P^2, \Phi_{16})$} &
{\scriptsize $\rho_1\rho_3(\rho_1+\rho_2),
\rho_1(\rho_1+\rho_2)(\rho_1+\rho_3),
\rho_1(\rho_1+\rho_2)(\rho_2+\rho_3),
 $}\\ & {\scriptsize
 $\rho_1(\rho_1+\rho_2)(\rho_1+\rho_2+\rho_3),
 \rho_3(\rho_1+\rho_2)(\rho_1+\rho_3), (\rho_1+\rho_2)(\rho_2+\rho_3)(\rho_1+\rho_2+\rho_3)$}\\
 \hline {\scriptsize $(S^1\times {\Bbb R}P^2, \Phi_{17})$} &
{\scriptsize $\rho_1\rho_3(\rho_1+\rho_2+\rho_3),
\rho_1(\rho_1+\rho_3)(\rho_1+\rho_2+\rho_3),(\rho_1+\rho_2)(\rho_1+\rho_3)(\rho_1+\rho_2+\rho_3),
 $}\\ & {\scriptsize
 $\rho_3(\rho_1+\rho_3)(\rho_1+\rho_2+\rho_3),
 (\rho_1+\rho_2)(\rho_2+\rho_3)(\rho_1+\rho_2+\rho_3), (\rho_1+\rho_3)(\rho_2+\rho_3)(\rho_1+\rho_2+\rho_3)$}\\
 \hline {\scriptsize $(S^1\times {\Bbb R}P^2, \Phi_{18})$} &
{\scriptsize $\rho_2\rho_3(\rho_1+\rho_3),
\rho_2(\rho_1+\rho_3)(\rho_2+\rho_3),
\rho_3(\rho_1+\rho_2)(\rho_1+\rho_3),
 $}\\ & {\scriptsize
 $\rho_3(\rho_1+\rho_3)(\rho_1+\rho_2+\rho_3),
 \rho_3(\rho_1+\rho_3)(\rho_2+\rho_3), (\rho_1+\rho_2)(\rho_1+\rho_3)(\rho_1+\rho_2+\rho_3)$}\\
 \hline {\scriptsize $(S^1\times {\Bbb R}P^2, \Phi_{19})$} &
{\scriptsize $\rho_2\rho_3(\rho_1+\rho_2),
\rho_2(\rho_1+\rho_2)(\rho_1+\rho_3),
\rho_2(\rho_1+\rho_2)(\rho_2+\rho_3),
 $}\\ & {\scriptsize
 $\rho_2(\rho_1+\rho_2)(\rho_1+\rho_2+\rho_3),
 \rho_3(\rho_1+\rho_2)(\rho_2+\rho_3), (\rho_1+\rho_2)(\rho_1+\rho_3)(\rho_1+\rho_2+\rho_3)$}\\
 \hline {\scriptsize $(S^1\times {\Bbb R}P^2, \Phi_{20})$} &
{\scriptsize $\rho_2\rho_3(\rho_1+\rho_2+\rho_3),
\rho_2(\rho_2+\rho_3)(\rho_1+\rho_2+\rho_3),(\rho_1+\rho_2)(\rho_1+\rho_3)(\rho_1+\rho_2+\rho_3),
 $}\\ & {\scriptsize
 $\rho_3(\rho_2+\rho_3)(\rho_1+\rho_2+\rho_3),
 (\rho_1+\rho_2)(\rho_2+\rho_3)(\rho_1+\rho_2+\rho_3), (\rho_1+\rho_3)(\rho_2+\rho_3)(\rho_1+\rho_2+\rho_3)$}\\
\hline
\end{tabular}

\vskip .2cm
 Now let $\beta\in
{\frak M}_3$ be an essential generator. By Lemma~\ref{l1}, one has
known that $\vert{\mathcal{N}}_\beta\vert\leq 14$.

\vskip .2cm

\noindent {\bf Claim 1.} {\em $\vert{\mathcal{N}}_\beta\vert$ must
be less than 12.}

\begin{proof} If
$\vert{\mathcal{N}}_\beta\vert=14$, then for each $i$
$(i=0,1,...,6)$, there must be two monomials $\delta^{(i)}_1,
\delta^{(i)}_2$ in $\mathcal{N}_{({\Bbb R}P^3, T_i)}$ such that
$\delta^{(i)}_1$ and $\delta^{(i)}_2$
 contain in ${\mathcal{N}}_\beta$. By Lemma~\ref{l6} and
 Remark~\ref{re}, an easy argument shows that there must be some  $(S^1\times {\Bbb R}P^2, \Phi_j)$ such that
 ${\mathcal{N}}_{(S^1\times {\Bbb
 R}P^2, \Phi_j)}\subset {\mathcal{N}}_\beta$. Then $8=\vert{\mathcal{N}}_{\beta+[(S^1\times {\Bbb
 R}P^2, \Phi_j)]}\vert<\vert{\mathcal{N}}_{\beta}\vert=14$. However, this is a
 contradiction since $\beta$ is an essential generator. Thus,
 $\vert{\mathcal{N}}_\beta\vert=14$ is impossible.

 \vskip .2cm

 If $\vert{\mathcal{N}}_\beta\vert=12$, since each $\mathcal{N}_{({\Bbb R}P^3, T_i)}$ contains at most two monomials in
 ${\mathcal{N}}_\beta$, then ${\mathcal{N}}_\beta=\{\delta_1,$ $\delta_2, ... , \delta_{11}, \delta_{12}\}$
 has the following two possible cases: (i) all elements of ${\mathcal{N}}_\beta$
 may be divided into six pairs $\{\delta_1,\delta_2\}$, $... , \{\delta_{11}, \delta_{12}\}$  that are  distributed in
 six different $\mathcal{N}_{({\Bbb R}P^3, T_i)}$ respectively; (ii) ${\mathcal{N}}_\beta$ may be divided into seven
 parts $\{\delta_1,\delta_2\}$, $... , \{\delta_{9}, \delta_{10}\}$,
 $\{\delta_{11}\}$, $\{\delta_{12}\}$ such that these seven parts are
 just distributed in $\mathcal{N}_{({\Bbb R}P^3, T_0)}, ..., \mathcal{N}_{({\Bbb R}P^3, T_6)}$, respectively.
 The similar argument as above also shows that there must be some  $(S^1\times {\Bbb R}P^2, \Phi_j)$ such that
 for the case (i), at least five elements of
 ${\mathcal{N}}_{(S^1\times {\Bbb
 R}P^2, \Phi_j)}$ contain in ${\mathcal{N}}_\beta$, and for the case (ii), at least four elements of
 ${\mathcal{N}}_{(S^1\times {\Bbb
 R}P^2, \Phi_j)}$ contain in ${\mathcal{N}}_\beta$. Then $\vert{\mathcal{N}}_{\beta+[(S^1\times {\Bbb
 R}P^2, \Phi_j)]}\vert\leq 10<\vert{\mathcal{N}}_{\beta}\vert=12$.
 This contradicts that $\beta$ is an essential generator.
 Thus,
 $\vert{\mathcal{N}}_\beta\vert=12$ cannot occur.
 \end{proof}

Let $(M, \phi)$ be a representative of $\beta$ such that
$\mathcal{N}_M$ is prime, and let $(\Gamma_M,\alpha)$ be the moment
graph of $(M, \phi)$.

\vskip .2cm \noindent {\bf Claim 2.} {\em $\Gamma_M$ is connected.}

 \begin{proof}
 Suppose that $\Gamma_M$ is disconnected. Let
$\Gamma^{'}$ be a connected component of $\Gamma_M$. Then the
restriction $\alpha|_{\Gamma^{'}}$ is still an axial function of
$\Gamma^{'}$. By Claim 1, one has $\vert{\mathcal{N}}_M\vert\leq
 10$ so the number of vertices of $\Gamma_M$ is less than or equal to 10. If $|V_{\Gamma^{'}}|=2$, then obviously
$\alpha(E_{p_1})=\alpha(E_{p_2})$ for $p_1,p_2\in V_{\Gamma^{'}}$,
but this is impossible since $\mathcal{N}_M$ is prime. If
$|V_{\Gamma^{'}}|=4$, then $\Gamma^{'}$ must be the 1-skeleton of a
3-simplex, and thus $(\Gamma^{'}, \alpha|_{\Gamma^{'}})$ is the
moment graph of some $({\Bbb R}P^3, T_i)$. Further, the disjoint
union of $(M, \phi)$ and $({\Bbb R}P^3, T_i)$ forms a $({\Bbb
Z}_2)^3$-action with at most six fixed points. This contradicts to
the assumption that $\beta$ is an essential generator. If
$|V_{\Gamma^{'}}|=6$, since the number of vertices of $\Gamma_M$ is
less than or equal to 10, $\Gamma_M$ must have another connected
component with 2 or 4 vertices, so that the problem is reduced to
the case $|V_{\Gamma^{'}}|=2$ or 4. This completes the proof.
 \end{proof}

\vskip .2cm By Proposition~\ref{gr} and Claim 2, the 2-skeletal
expansion $N$ of $(\Gamma_M,\alpha)$ is a connected closed surface.
 By $F_{\Gamma_M}$ one denotes the
set of all
 2-nests in $\mathcal{K}_{(\Gamma_M,\alpha)}$. Then one has the
 formula
 \begin{eqnarray}\chi(N)=\vert V_{\Gamma_M}\vert-\vert E_{\Gamma_M}\vert+
 \vert F_{\Gamma_M}\vert
 \end{eqnarray}
 where $\chi(N)$ is the Euler characteristic of $N$.
 Note that $\vert
 V_{\Gamma_M}\vert=\vert{\mathcal{N}}_\beta\vert$ and $
 3\vert V_{\Gamma_M}\vert=2\vert E_{\Gamma_M}\vert$.

 \vskip .2cm

 \noindent {\bf Claim 3.} {\em The 2-skeletal expansion $N$ of
$(\Gamma_M,\alpha)$ is a
 sphere of dimension 2.}

 \begin{proof}
It suffices to show that the Euler characteristic $\chi(N)$ is 2. By
Claim 1, one has $\vert{\mathcal{N}}_M\vert\leq
 10$ so one needs to consider the cases of
 $\vert{\mathcal{N}}_M\vert=4,6,8,10$.

 \vskip .2cm

 When $\vert{\mathcal{N}}_M\vert=4$, if $\chi(N)$ is not  2, then from (6.1)
 one has that $\vert
 F_{\Gamma_M}\vert\leq 3$, so all 2-nests in $(\Gamma_M, \alpha)$ correspond to at most three
 nonzero elements in $\H({\Bbb Z}_2, ({\Bbb Z}_2)^3)$. However,
 any three nonzero elements in $\H({\Bbb Z}_2, ({\Bbb Z}_2)^3)$ cannot produce four different bases
 of  $\H({\Bbb Z}_2, ({\Bbb Z}_2)^3)$. Thus, $\chi(N)$ must be 2.

 \vskip .2cm

 When $\vert{\mathcal{N}}_M\vert=6$, since any four
  nonzero elements in $\H({\Bbb Z}_2, ({\Bbb Z}_2)^3)$ cannot produce six different bases
 of  $\H({\Bbb Z}_2, ({\Bbb Z}_2)^3)$, one has that $\vert
 F_{\Gamma_M}\vert$ must be 5 so $\chi(N)$ is 2.

 \vskip .2cm

 When
 $\vert{\mathcal{N}}_M\vert=8$, if $N$ is
 not a
 sphere of dimension 2, then the above argument makes sure that $\vert
 F_{\Gamma_M}\vert$ must be 5, and the dual map $\eta$ of $\alpha$ maps five 2-nests of $\mathcal{K}_{(\Gamma_M,
 \alpha)}$  to five different nonzero elements of $\H({\Bbb Z}_2, ({\Bbb Z}_2)^3)$, respectively.
 An easy argument shows that any  five   nonzero elements in $\H({\Bbb Z}_2, ({\Bbb Z}_2)^3)$
 can be translated into five given nonzero elements by applying an
 automorphism of $\H({\Bbb Z}_2, ({\Bbb Z}_2)^3)$. Thus we may
 choose five special elements $\rho_1^*, \rho_2^*, \rho_3^*,
 \rho_1^*+\rho_2^*, \rho_1^*+\rho_3^*$ as being the images of
 $\eta$ on five 2-nests of $\mathcal{K}_{(\Gamma_M,
 \alpha)}$, where $\{\rho_1^*, \rho_2^*, \rho_3^*\}$ is the
 standard basis of $\H({\Bbb Z}_2, ({\Bbb Z}_2)^3)$, which
 corresponds to the standard basis $\{\rho_1, \rho_2, \rho_3\}$ of
 $\H(({\Bbb Z}_2)^3, {\Bbb Z}_2)$. Then from these five chosen elements, one may produce just 8 bases of
 $\H({\Bbb Z}_2, ({\Bbb Z}_2)^3)$ as follows:
 $$\{\rho_1^*,\rho_2^*,\rho_3^*\},
 \{\rho_1^*,\rho_2^*,\rho_1^*+\rho_3^*\},
 \{\rho_1^*,\rho_3^*,\rho_1^*+\rho_2^*\},
 \{\rho_1^*,\rho_1^*+\rho_2^*,\rho_1^*+\rho_3^*\},$$
 $$\{\rho_2^*,\rho_3^*,\rho_1^*+\rho_2^*\},
 \{\rho_2^*,\rho_3^*,\rho_1^*+\rho_3^*\},
 \{\rho_2^*,\rho_1^*+\rho_2^*,\rho_1^*+\rho_3^*\},
 \{\rho_3^*,\rho_1^*+\rho_2^*,\rho_1^*+\rho_3^*\}.$$
So, $\mathcal{N}_M$ consists of 8 monomials $\rho_1\rho_2\rho_3$,
$\rho_2\rho_3(\rho_1+\rho_3)$, $\rho_2\rho_3(\rho_1+\rho_2)$,
$\rho_2\rho_3(\rho_1+\rho_2+\rho_3)$,
$\rho_1\rho_3(\rho_1+\rho_2)$, $\rho_1\rho_2(\rho_1+\rho_3)$,
$\rho_3(\rho_1+\rho_3)(\rho_1+\rho_2+\rho_3)$,
$\rho_2(\rho_1+\rho_2)(\rho_1+\rho_2+\rho_3)$. Further, we see
from Table I that $\mathcal{N}_{({\Bbb R}P^3, T_3)}\subset
\mathcal{N}_M$, so $|\mathcal{N}_{\beta+[({\Bbb R}P^3, T_3)]}|<8$.
This means that $\beta$ is not an essential generator, which gives
a contradiction. Thus, when
 $\vert{\mathcal{N}}_M\vert=8$,
 $\vert
 F_{\Gamma_M}\vert$ must be 6 so $\chi(N)$ is still 2.

 \vskip .2cm

 When $\vert{\mathcal{N}}_M\vert=10$, suppose that  $\chi(N)$ is not  2.
 As shown above, any
 five
  nonzero elements in $\H({\Bbb Z}_2, ({\Bbb Z}_2)^3)$ cannot produce ten different bases
 of  $\H({\Bbb Z}_2, ({\Bbb Z}_2)^3)$, and thus the only
 possibility of  $\vert
 F_{\Gamma_M}\vert$ is 6. Further, one has from (6.1) that $\chi(N)$
 must be 1. To ensure that $\vert{\mathcal{N}}_M\vert=10$, six 2-nests in $\mathcal{K}_{(\Gamma_M, \alpha)}$ must
 then
 correspond to six different nonzero elements in $\H({\Bbb Z}_2, ({\Bbb Z}_2)^3)$ by the dual map $\eta$.
 It is easy to check that
  any six different
 nonzero elements in $\H({\Bbb Z}_2, ({\Bbb Z}_2)^3)$ can still be translated into the given six
 different nonzero
 elements by an automorphism of $\H({\Bbb Z}_2, ({\Bbb Z}_2)^3)$.  Thus, as in the argument of the case
 $\vert{\mathcal{N}}_M\vert=8$,  it needs to merely
 consider six special nonzero elements $\H({\Bbb Z}_2, ({\Bbb Z}_2)^3)$. Take six nonzero elements
  $\rho_1^*, \rho^*_2, \rho^*_3, \rho^*_1+\rho^*_2,
 \rho^*_1+\rho^*_3, \rho^*_1+\rho^*_2+\rho^*_3$ in $\H({\Bbb Z}_2, ({\Bbb Z}_2)^3)$ such that they are the images
 of $\eta$ on six 2-nests, one then may
 produce 16 different bases of $\H({\Bbb Z}_2, ({\Bbb Z}_2)^3)$ as follows:
 $\{\rho_1^*, \rho_2^*, \rho_3^*\}, \{\rho_2^*, \rho_3^*, \rho_1^*+\rho_2^*+\rho_3^*\}, \{\rho_1^*, \rho_3^*,
 \rho_1^*+\rho_2^*+\rho_3^*\},
 \{\rho_1^*, \rho_2^*, \rho_1^*+\rho_2^*+\rho_3^*\},$
 $\{\rho_3^*, \rho_1^*+\rho_2^*, \rho_1^*+\rho_3^*\}, \{\rho_2^*, \rho_3^*, \rho_1^*+\rho_3^*\}, \{\rho_2^*,
 \rho_1^*+\rho_2^*, \rho_1^*+\rho_3^*\},
 \{\rho_2^*, \rho_3^*, \rho_1^*+\rho_2^*\},$
$\{\rho_1^*+\rho_2^*, \rho_1^*+\rho_3^*,
\rho_1^*+\rho_2^*+\rho_3^*\}, \{\rho_1^*, \rho_1^*+\rho_3^*,
\rho_1^*+\rho_2^*+\rho_3^*\}, \{\rho_1^*, \rho_1^*+\rho_2^*,
\rho_1^*+\rho_2^*+\rho_3^*\},
 \{\rho_1^*, \rho_1^*+\rho_2^*, \rho_1^*+\rho_3^*\},$
 $\{\rho_1^*, \rho_3^*, \rho_1^*+\rho_2^*\}, \{\rho_1^*, \rho_2^*, \rho_1^*+\rho_3^*\},
 \{\rho_3^*, \rho_1^*+\rho_3^*, \rho_1^*+\rho_2^*+\rho_3^*\},
 \{\rho_2^*, \rho_1^*+\rho_2^*, \rho_1^*+\rho_2^*+\rho_3^*\}.$
 These 16 bases are dual to 16 bases in $\H(({\Bbb Z}_2)^3, {\Bbb
 Z}_2)$, which give the following 16   monomials.
$$\rho_1\rho_2\rho_3, \rho_1(\rho_1+\rho_2)(\rho_1+\rho_3), \rho_2(\rho_1+\rho_2)(\rho_2+\rho_3),
\rho_1\rho_2(\rho_1+\rho_2+\rho_3),$$
$$\rho_2(\rho_1+\rho_2)(\rho_1+\rho_2+\rho_3), \rho_1\rho_2(\rho_1+\rho_3), \rho_3(\rho_1+\rho_3)(\rho_1+\rho_2+\rho_3),
\rho_1\rho_3(\rho_1+\rho_2),$$
$$(\rho_1+\rho_2)(\rho_1+\rho_3)(\rho_1+\rho_2+\rho_3), \rho_2(\rho_1+\rho_3)(\rho_2+\rho_3),
\rho_3(\rho_1+\rho_2)(\rho_2+\rho_3),
\rho_2\rho_3(\rho_1+\rho_2+\rho_3),$$
$$\rho_2\rho_3(\rho_1+\rho_2), \rho_2\rho_3(\rho_1+\rho_3), \rho_2(\rho_1+\rho_2)(\rho_1+\rho_3),
\rho_3(\rho_1+\rho_2)(\rho_1+\rho_3).$$ One sees  that the first
row above is just $\mathcal{N}_{({\Bbb R}P^3, T_0)}$, the second
row is $\mathcal{N}_{({\Bbb R}P^3, T_3)}$, and the third row is
$\mathcal{N}_{({\Bbb R}P^3, T_6)}$, but
$\rho_2\rho_3(\rho_1+\rho_2), \rho_2\rho_3(\rho_1+\rho_3),
\rho_2(\rho_1+\rho_2)(\rho_1+\rho_3),
\rho_3(\rho_1+\rho_2)(\rho_1+\rho_3)$ belong to
$\mathcal{N}_{({\Bbb R}P^3, T_1)}$, $\mathcal{N}_{({\Bbb R}P^3,
T_2)}$, $\mathcal{N}_{({\Bbb R}P^3, T_4)}$, $\mathcal{N}_{({\Bbb
R}P^3, T_5)}$, respectively. Then
 $\mathcal{N}_M$ must contain $\rho_2\rho_3(\rho_1+\rho_2), \rho_2\rho_3(\rho_1+\rho_3),
\rho_2(\rho_1+\rho_2)(\rho_1+\rho_3),
\rho_3(\rho_1+\rho_2)(\rho_1+\rho_3)$, and $|\mathcal{N}_M\cap
\mathcal{N}_{({\Bbb R}P^3, T_i)}|=2$ for $i=0,3,6$.

\vskip .2cm

Now choose any two $\gamma_1, \gamma_2$ of
$\rho_2\rho_3(\rho_1+\rho_2), \rho_2\rho_3(\rho_1+\rho_3),
\rho_2(\rho_1+\rho_2)(\rho_1+\rho_3),
\rho_3(\rho_1+\rho_2)(\rho_1+\rho_3)$, it is easy to show that there
is always
 one  $(S^1\times{\Bbb R}P^2, \Phi_j)$ such that $\mathcal{N}_{(S^1\times{\Bbb R}P^2, \Phi_j)}$ contains $\gamma_1, \gamma_2$.
Without loss of generality, we may let
$\gamma_1=\rho_2\rho_3(\rho_1+\rho_2)$ and
$\gamma_2=\rho_2\rho_3(\rho_1+\rho_3)$. Then one has that
$\mathcal{N}_{(S^1\times{\Bbb R}P^2,
\Phi_j)}=\{\rho_2\rho_3(\rho_1+\rho_2),
\rho_2\rho_3(\rho_1+\rho_3),\rho_1\rho_2\rho_3,
\rho_1\rho_2(\rho_1+\rho_3),\rho_2\rho_3(\rho_1+\rho_2+\rho_3),\rho_2(\rho_1+\rho_2)(\rho_1+\rho_2+\rho_3)
\}$ with $\rho_1\rho_2\rho_3\in \mathcal{N}_{({\Bbb R}P^3, T_0)}$,
$\rho_1\rho_2(\rho_1+\rho_3),
\rho_2(\rho_1+\rho_2)(\rho_1+\rho_2+\rho_3)\in \mathcal{N}_{({\Bbb
R}P^3, T_3)}$, $\rho_2\rho_3(\rho_1+\rho_2+\rho_3)\in
\mathcal{N}_{({\Bbb R}P^3, T_6)}$. If $\mathcal{N}_M$ contains at
least two of $\rho_1\rho_2\rho_3,
\rho_1\rho_2(\rho_1+\rho_3),\rho_2\rho_3(\rho_1+\rho_2+\rho_3),\rho_2(\rho_1+\rho_2)(\rho_1+\rho_2+\rho_3)$,
form the disjoint union $(M,\phi)\sqcup(S^1\times{\Bbb R}P^2,
\Phi_j)$, then
$$|\mathcal{N}_{\beta+[(S^1\times{\Bbb R}P^2, \Phi_j)]}|<10,$$
which contradicts that $\beta$ is an essential generator. Thus,
this case cannot occur. If $\mathcal{N}_M$ contains only one (say
$\omega$) of $\rho_1\rho_2\rho_3,
\rho_1\rho_2(\rho_1+\rho_3),\rho_2\rho_3(\rho_1+\rho_2+\rho_3),\rho_2(\rho_1+\rho_2)(\rho_1+\rho_2+\rho_3)$,
form the union $(M,\phi)\sqcup(S^1\times{\Bbb R}P^2,
\Phi_j)\sqcup({\Bbb R}P^3, T_l)$  where $({\Bbb R}P^3, T_l)$ is
some one of $({\Bbb R}P^3, T_0), ({\Bbb R}P^3, T_3), ({\Bbb R}P^3,
T_6)$ such that $\omega\not\in \mathcal{N}_{({\Bbb R}P^3, T_l)}$,
then
$$|\mathcal{N}_{\beta+[(S^1\times{\Bbb R}P^2, \Phi_j)]+[({\Bbb R}P^3, T_l)]}|<10,$$
which leads to a contradiction (note that
$|\mathcal{N}_{[(S^1\times{\Bbb R}P^2, \Phi_j)]+[({\Bbb R}P^3,
T_l)]}|<10$). Finally, if $\mathcal{N}_M$ does not contain  any one
of $\rho_1\rho_2\rho_3,
\rho_1\rho_2(\rho_1+\rho_3),\rho_2\rho_3(\rho_1+\rho_2+\rho_3),\rho_2(\rho_1+\rho_2)(\rho_1+\rho_2+\rho_3)$,
consider the disjoint union of $(M,\phi)\sqcup(S^1\times{\Bbb R}P^2,
\Phi_j)$ with  $({\Bbb R}P^3, T_3)$, then a contradiction still
occurs (i.e., $\beta$ is not an essential generator). This is
impossible. Therefore, $\chi(N)$ must be  2.

\vskip .2cm

Combining the above arguments, we complete the proof.
 \end{proof}
 \begin{lem} \label{l7}
Let $\beta\in {\frak M}_3$. If $\beta$ is an essential generator,
then $\vert{\mathcal{N}}_\beta\vert\leq 6$.
\end{lem}
\begin{proof}
By Claim 1 it suffices to show that $\vert{\mathcal{N}}_\beta\vert$
is not equal to $8$ and $10$. One knows by Claim 3 that the
2-skeletal expansion $N$ is a sphere of dimension 2, so  $\Gamma_M$
is planar and in particular, it is the 1-skeleton  of a simple
convex 3-polytope $P^3$. In this case, $M$ is a small cover over
$P^3$, so the axial function $\alpha$ on $\Gamma_M$ is dual to the
characteristic function $\lambda$ on $P^3$.

The argument proceeds as follows.

\vskip .2cm

{\em Case $(i)$}: $\vert{\mathcal{N}}_\beta\vert=8$.

\vskip .2cm If $\vert{\mathcal{N}}_\beta\vert=8$, then $\Gamma_M$ is
the 1-skeleton of a simple convex polytope with 8 vertices. From
\cite{g} one knows that there are only two different combinatorial
types of simple 3-polytopes with eight vertices, as shown in
Figure~\ref{a2}.
\begin{figure}[h]
    \input{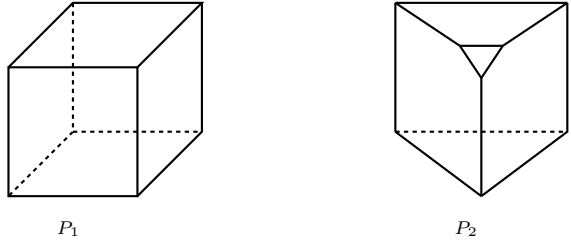}\centering
    \caption[]{Two simple 3-polytopes with eight
vertices }\label{a2}
\end{figure}
If $\Gamma_M$ is the 1-skeleton of 3-dimensional cube $P_1$, then
it is easy to check that $P_1$ does not admit any characteristic
function of mapping six 2-faces into  six different nonzero
elements in $\H({\Bbb Z}_2, ({\Bbb Z}_2)^3)$, but this is
impossible. Thus, $\Gamma_M$ cannot be the 1-skeleton of
 $P_1$. If $\Gamma_M$ is the 1-skeleton of
 $P_2$, taking  a  triangular facet $F$ of $P_2$, then, up to automorphisms of $\H({\Bbb Z}_2, ({\Bbb
 Z}_2)^3)$, it is easy to see that the characteristic
 function $\lambda$ on $P_2$ maps $F$ with its 3 adjacent 2-faces into either $\rho_1^*, \rho_2^*, \rho_3^*,
\rho_1^*+\rho_2^*+\rho_3^*$ or $\rho_1^*, \rho_2^*, \rho_3^*,
\rho_1^*+\rho_2^*$. When $\lambda$ maps $F$ with its 3 adjacent
2-faces into  $\rho_1^*, \rho_2^*, \rho_3^*,
\rho_1^*+\rho_2^*+\rho_3^*$, obviously there must be some $({\Bbb
R}P^3, T_i)$ such that $|\mathcal{N}_{\beta+[({\Bbb R}P^3,
T_i)]}|=6<8$. This contradicts the fact that $\beta$ is an
essential generator, and thus this case cannot occur. When
$\lambda$ maps $F$ with its 3 adjacent 2-faces into  $\rho_1^*,
\rho_2^*, \rho_3^*, \rho_1^*+\rho_2^*$, it is easy to check that
there must be some $(S^1\times{\Bbb R}P^2, \Phi_j)$ such that
$|\mathcal{N}_{\beta+[(S^1\times{\Bbb R}P^2, \Phi_j)]}|=6<8$. This
also is impossible, so $\Gamma_M$ cannot be  the 1-skeleton of
 $P_2$. Thus, if $\beta$ is an essential generator, then
 $\vert{\mathcal{N}}_\beta\vert=8$ is impossible.

\vskip .2cm

{\em Case $(ii)$}: $\vert{\mathcal{N}}_\beta\vert=10$.

\vskip .2cm If $\vert{\mathcal{N}}_\beta\vert=10$, then $\Gamma_M$
is the 1-skeleton of a simple convex polytope with 10 vertices. From
\cite{g} one knows that there are only five different combinatorial
types of simple 3-polytopes with ten vertices, as shown in
Figures~\ref{a3} and~\ref{a4}.
\begin{figure}[h]
    \input{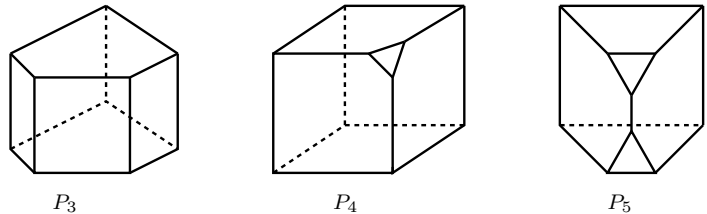}\centering
    \caption[]{Simple 3-polytopes with ten
vertices }\label{a3}
\end{figure}
 \begin{figure}[h]
    \input{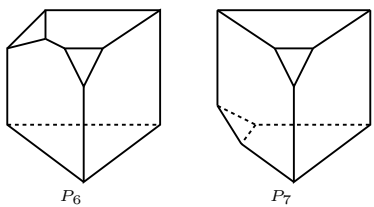}\centering
    \caption[]{ Simple 3-polytopes with ten
vertices }\label{a4}
\end{figure}
An easy argument shows that $\Gamma_M$ cannot be the 1-skeleton of
$P_3$. Since each of $P_4,P_5,P_6,P_7$ has at least one triangular
facet, similarly to the proof of case (i), one may prove that
$\Gamma_M$ cannot be the 1-skeleton of $P_4,P_5,P_6,P_7$,
respectively. Therefore,   $\vert{\mathcal{N}}_\beta\vert=10$
 is impossible, too.

 \vskip .2cm

 Combining the above arguments, one completes the proof.
\end{proof}

Together with Lemma~\ref{l5}, Lemma~\ref{l7} and Remark~\ref{re},
we complete the proof of Proposition~\ref{p}.

\begin{thm}\label{dim}
As a vector space over ${\Bbb Z}_2$,  ${\frak M}_3$ has dimension
13, and it is generated by $({\Bbb R}P^3, T_0), ({\Bbb R}P^3,
T_1),...,({\Bbb R}P^3, T_6), (S^1\times {\Bbb R}P^2, \Phi_0),
(S^1\times {\Bbb R}P^2, \Phi_1), ..., (S^1\times {\Bbb R}P^2,
\Phi_4), (S^1\times {\Bbb R}P^2, \Phi_6)$.
\end{thm}
\begin{proof} By Propositions~\ref{l2} and~\ref{p}, any element of ${\frak
M}_3$ is a linear combination of 28 small covers $({\Bbb R}P^3,
T_i), i=0,1,...,6$, and $(S^1\times {\Bbb R}P^2, \Phi_j),
j=0,1,...,20$. Thus, in order to calculate the dimension of
${\frak M}_3$, one needs to determine a maximal linearly
independent set of  the above 28 small covers. Let
$$\sum_{i=0}^6l_i[({\Bbb
R}P^3, T_i)]+\sum_{j=0}^{20}k_j[(S^1\times{\Bbb R}P^2),
\Phi_j)]=0$$ where $l_i,k_j\in {\Bbb Z}_2$. Using
Stong-homomorphism  $\delta_3$ in Theorem~\ref{s}, one then has
that \begin{equation} \label{e}\sum_{i=0}^6l_i\delta_3([({\Bbb
R}P^3, T_i)])+\sum_{j=0}^{20}k_j\delta_3([(S^1\times{\Bbb R}P^2),
\Phi_j)])=0. \end{equation} Since $\H(({\Bbb Z}_2)^3,{\Bbb Z}_2)$
gives 28 different bases, from (\ref{e}) and Tables I and II,  one
obtains a equation system  formed by 28 equations, such that the
coefficient matrix of this equation system is

 {\tiny
\begin{eqnarray*} \setcounter{MaxMatrixCols}{28} A=
\begin{pmatrix}
1 & 0 &0 &0 &0  &0&0  &1 & 1& 1&1 &0 & 0& 1& 0& 0& 1& 0 &0& 0& 0&0 &0 &0 &0 &0 &0 &0 \\
  1 &0&0&0&0&0&0&1&0&0&0&0&1&0&0&0&0&0&0&0&1&0&0&1&0&0&0&0 \\
   1 & 0 &0 &0 &0  &0 &0&0 &1 &0&0&0&0&0&0&1&0&0&0&1&0&0&0&0&0&0&1&0 \\
  1 &0&0&0&0&0&0&0&0&1&0&0&0&0&0&0&0&0&1&0&0&0&1&0&0&1&0&0 \\
  0&1&0&0&0&0&0&0&0&0&1&1&1&0&1&0&0&0&0&0&0&1&0&0&0&0&1&0\\
  0&1&0&0&0&0&0&0&1&0&1&0&0&0&0&0&1&0&0&0&1&0&0&0&0&0&0&0\\
  0&1&0&0&0&0&0&0&0&0&0&1&0&0&0&0&0&1&0&0&0&0&0&0&1&1&0&0\\
  0&1&0&0&0&0&0&0&0&1&0&0&1&1&0&0&0&0&0&0&0&0&0&1&0&0&0&0\\
  0&0&1&0&0&0&0&0&0&0&0&1&0&1&1&1&0&0&0&0&0&1&0&1&0&0&0&0\\
  0&0&1&0&0&0&0&1&0&0&0&0&0&1&0&0&1&0&0&1&0&0&0&0&0&0&0&0\\
  0&0&1&0&0&0&0&0&0&0&0&0&0&0&1&0&0&1&0&0&0&0&1&0&0&0&0&1\\
  0&0&1&0&0&0&0&0&0&1&1&0&0&0&0&1&0&0&0&0&0&0&0&0&0&0&1&0\\
  0&0&0&1&0&0&0&0&0&0&0&1&0&0&0&0&1&1&1&0&1&0&0&0&1&0&0&0\\
  0&0&0&1&0&0&0&1&0&0&0&0&0&1&0&0&1&0&0&0&0&0&1&0&0&0&0&0\\
  0&0&0&1&0&0&0&0&0&0&0&0&0&0&1&0&0&1&0&1&0&0&0&0&0&0&0&1\\
  0&0&0&1&0&0&0&0&1&0&1&0&0&0&0&0&0&0&1&0&0&0&0&0&0&1&0&0\\
  0&0&0&0&1&0&0&1&0&0&0&0&0&0&0&1&0&0&0&1&1&1&0&1&0&0&0&0\\
  0&0&0&0&1&0&0&0&0&0&0&0&0&0&0&0&0&0&1&1&0&0&0&0&1&0&0&1\\
  0&0&0&0&1&0&0&0&0&0&0&1&0&0&1&0&1&0&0&0&1&0&0&0&0&0&0&0\\
  0&0&0&0&1&0&0&0&1&0&0&0&1&0&0&0&0&0&0&0&0&1&0&0&0&0&1&0\\
  0&0&0&0&0&1&0&0&0&0&0&0&0&0&0&1&0&0&0&0&0&1&1&1&1&0&0&1\\
  0&0&0&0&0&1&0&1&0&0&0&0&0&0&0&0&0&0&1&0&1&0&1&0&0&0&0&0\\
  0&0&0&0&0&1&0&0&0&1&0&0&1&0&0&0&0&0&0&0&0&0&0&1&0&1&0&0\\
  0&0&0&0&0&1&0&0&0&0&0&1&0&1&0&0&0&1&0&0&0&0&0&0&1&0&0&0\\
  0&0&0&0&0&0&1&0&0&0&0&0&1&0&0&0&0&0&0&0&0&1&0&0&1&1&1&1\\
  0&0&0&0&0&0&1&0&1&0&0&0&0&0&0&0&0&0&1&1&0&0&0&0&0&1&0&0\\
  0&0&0&0&0&0&1&0&0&1&0&0&0&0&0&1&0&0&0&0&0&0&1&0&0&0&1&0\\
  0&0&0&0&0&0&1&0&0&0&1&0&0&0&1&0&0&1&0&0&0&0&0&0&0&0&0&1
 \end{pmatrix}.
    \end{eqnarray*}}
    By doing elementary row operations,  $A$ is changed into
 {\tiny   \begin{eqnarray*}
\setcounter{MaxMatrixCols}{28} A^{'}=
\begin{pmatrix}
1&0&0&0&0&0&0  &0&0&1&0&0&0&0&0&0&0&0&1&0&0&0&1&0&0&1&0&0 \\
0&1&0&0&0&0&0  &0&0&0&0&1&0&0&0&0&0&1&0&0&0&0&0&0&1&1&0&0 \\
0&0&1&0&0&0&0  &0&0&0&0&0&0&0&1&0&0&1&0&0&0&0&1&0&0&0&0&1 \\
0&0&0&1&0&0&0  &0&0&0&0&0&0&0&1&0&0&1&0&1&0&0&0&0&0&0&0&1 \\
0&0&0&0&1&0&0  &0&0&0&0&0&0&0&0&0&0&0&1&1&0&0&0&0&1&0&0&1\\
0&0&0&0&0&1&0  &0&0&0&0&0&0&0&0&1&0&0&0&0&0&1&1&1&1&0&0&1\\
0&0&0&0&0&0&1  &0&0&0&0&0&1&0&0&0&0&0&0&0&0&1&0&0&1&1&1&1\\
0&0&0&0&0&0&0  &1&0&0&0&0&0&0&0&1&0&0&1&0&1&1&0&1&1&0&0&1\\
0&0&0&0&0&0&0  &0&1&0&0&0&1&0&0&0&0&0&1&1&0&1&0&0&1&0&1&1\\
0&0&0&0&0&0&0  &0&0&1&0&0&1&0&0&1&0&0&0&0&0&1&1&0&1&1&0&1\\
0&0&0&0&0&0&0  &0&0&0&1&0&1&1&0&1&1&0&1&1&1&0&1&1&0&1&1&0\\
0&0&0&0&0&0&0  &0&0&0&0&1&0&0&1&0&1&0&1&1&1&0&0&0&1&0&0&1\\
0&0&0&0&0&0&0  &0&0&0&0&0&0&1&1&1&1&1&1&1&1&1&1&1&1&0&0&0\\
  0&0&0&0&0&0&0&0&0&0&0&0&0&0&0&0&0&0&0&0&0&0&0&0&0&0&0&0\\
  0&0&0&0&0&0&0&0&0&0&0&0&0&0&0&0&0&0&0&0&0&0&0&0&0&0&0&0\\
  0&0&0&0&0&0&0&0&0&0&0&0&0&0&0&0&0&0&0&0&0&0&0&0&0&0&0&0\\
  0&0&0&0&0&0&0&0&0&0&0&0&0&0&0&0&0&0&0&0&0&0&0&0&0&0&0&0\\
  0&0&0&0&0&0&0&0&0&0&0&0&0&0&0&0&0&0&0&0&0&0&0&0&0&0&0&0\\
  0&0&0&0&0&0&0&0&0&0&0&0&0&0&0&0&0&0&0&0&0&0&0&0&0&0&0&0\\
  0&0&0&0&0&0&0&0&0&0&0&0&0&0&0&0&0&0&0&0&0&0&0&0&0&0&0&0\\
  0&0&0&0&0&0&0&0&0&0&0&0&0&0&0&0&0&0&0&0&0&0&0&0&0&0&0&0\\
  0&0&0&0&0&0&0&0&0&0&0&0&0&0&0&0&0&0&0&0&0&0&0&0&0&0&0&0\\
  0&0&0&0&0&0&0&0&0&0&0&0&0&0&0&0&0&0&0&0&0&0&0&0&0&0&0&0\\
  0&0&0&0&0&0&0&0&0&0&0&0&0&0&0&0&0&0&0&0&0&0&0&0&0&0&0&0\\
  0&0&0&0&0&0&0&0&0&0&0&0&0&0&0&0&0&0&0&0&0&0&0&0&0&0&0&0\\
  0&0&0&0&0&0&0&0&0&0&0&0&0&0&0&0&0&0&0&0&0&0&0&0&0&0&0&0\\
  0&0&0&0&0&0&0&0&0&0&0&0&0&0&0&0&0&0&0&0&0&0&0&0&0&0&0&0\\
  0&0&0&0&0&0&0&0&0&0&0&0&0&0&0&0&0&0&0&0&0&0&0&0&0&0&0&0\\
 \end{pmatrix}.
    \end{eqnarray*}}
    so the rank of $A$ is 13, which is just the dimension of ${\frak
M}_3$. Theorem~\ref{dim} then follows from this.
\end{proof}

\section{Representatives of equivariant cobordism classes of
${\frak M}_3$}

Given two small covers $\pi_i: M_i^n\longrightarrow P_i^n, i=1,2$,
their equivariant connected sum along fixed points can be
proceeded as follows: Take a vertex $v_i$ from $P_i^n$ and let
$p_i$ be its preimage in $M_i$, $i=1,2$. With no loss one may
assume that $({\Bbb Z}_2)^n$-actions are equivalent in a
neighborhood of $p_i$ (actually, if necessary, one can change the
action by using an automorphism of $({\Bbb Z}_2)^n$). Then one can
perform the connected sum equivariantly near the fixed points
$p_1, p_2$. The result is a  2-torus manifold $M_1^n\sharp M_2^n$,
and its orbit space, $P_1^n\sharp P_2^n$, is given by removing a
small ball around $v_i$ from $P_i^n$ and gluing the results
together. As pointed out in \cite{dj}, generally $P_1^n\sharp
P_2^n$ is not canonically identified with a simple polytope but is
almost as good in that its boundary complex is dual to some PL
triangulation of $S^{n-1}$. However, it is easy to see that if
$n=3$, $P_1^n\sharp P_2^n$ is also a simple polytope, so
$M_1^n\sharp M_2^n$ is a small cover over $P_1^n\sharp P_2^n$.

\begin{lem} \label{l8} There exists a 3-dimensional small cover
$\pi: M^3\longrightarrow P^3$ such that $M$ is equivariantly
cobordant to a  2-torus 3-manifold $N$ with $\mathcal{N}_M$ prime
and $|\mathcal{N}_N|=28$.
\end{lem}

\begin{proof}
Consider two small covers $(S^1\times{\Bbb R}P^2, \Phi_0)$ and
$(S^1\times{\Bbb R}P^2, \Phi_1)$ over a prism $P^3$, one sees from
Table II that they have fixed points with the same representation
$\rho_1\rho_2\rho_3$. Then one can make an equivariant connected sum
along the fixed points with representation $\rho_1\rho_2\rho_3$,
such that $(S^1\times{\Bbb R}P^2, \Phi_0)\sharp(S^1\times{\Bbb
R}P^2, \Phi_1)$ is also a small cover over a simple 3-polytope with
10 vertices, and its tangent representation set is just equal to
$\mathcal{N}_{[(S^1\times{\Bbb R}P^2, \Phi_0)]+[(S^1\times{\Bbb
R}P^2, \Phi_1)]}$, consisting of $\rho_1\rho_2(\rho_2+\rho_3),
\rho_1\rho_3(\rho_2+\rho_3),\rho_1(\rho_1+\rho_2)(\rho_1+\rho_3),
\rho_1(\rho_1+\rho_2)(\rho_2+\rho_3),
\rho_1(\rho_1+\rho_3)(\rho_2+\rho_3), \rho_1\rho_2(\rho_1+\rho_3),
\rho_2\rho_3(\rho_1+\rho_3), \rho_2(\rho_1+\rho_2)(\rho_1+\rho_3),
\rho_2(\rho_1+\rho_2)(\rho_2+\rho_3),
\rho_2(\rho_1+\rho_3)(\rho_2+\rho_3)$. From Table I one sees the
following properties:

\vskip .2cm

(a)  For any $({\Bbb R}P^3, T_i)$, the intersection of
$\mathcal{N}_{[(S^1\times{\Bbb R}P^2, \Phi_0)]+[(S^1\times{\Bbb
R}P^2, \Phi_1)]}$ and $\mathcal{N}_{[({\Bbb R}P^3, T_i)]}$ is
always non-empty.

(b)   Two elements $\rho_1(\rho_1+\rho_2)(\rho_1+\rho_3),
\rho_2(\rho_1+\rho_2)(\rho_2+\rho_3)$ of
$\mathcal{N}_{[(S^1\times{\Bbb R}P^2, \Phi_0)]+[(S^1\times{\Bbb
R}P^2, \Phi_1)]}$ contain in $\mathcal{N}_{[({\Bbb R}P^3, T_0)]}$.

\vskip .2cm

Next, one preforms an equivariant connected sum of two copies of
$(S^1\times{\Bbb R}P^2, \Phi_0)\sharp(S^1\times{\Bbb R}P^2,
\Phi_1)$ along the fixed point with representation
$\rho_1(\rho_1+\rho_2)(\rho_1+\rho_3)$. Then the resulting $({\Bbb
Z}_2)^3$-manifold $M^{'}$ fixes 18 isolated points and  is also a
small cover over a simple polytope with 18 vertices. Obviously,
the representations at 18 fixed points of $M^{'}$ appear in pairs,
so $M^{'}$ bounds equivariantly. Since
$\mathcal{N}_{[(S^1\times{\Bbb R}P^2, \Phi_0)]+[(S^1\times{\Bbb
R}P^2, \Phi_1)]}\setminus
\{\rho_1(\rho_1+\rho_2)(\rho_1+\rho_3)\}\subset
\mathcal{N}_{M^{'}}$, by properties (a) and (b), one has that for
any $({\Bbb R}P^3, T_i)$, the intersection
$\mathcal{N}_{M^{'}}\cap \mathcal{N}_{[({\Bbb R}P^3, T_i)]}$ is
non-empty, so that $M^{'}$ can perform an equivariant connected
sum with each $({\Bbb R}P^3, T_i)$ along the fixed points with the
same representation. Let $M$ be the equivariant connected sum of
$M^{'}$ with all $({\Bbb R}P^3, T_i)$ in the above way. Then $M$
is just a desired small cover.
\end{proof}

\begin{thm}\label{small}
Any element $\beta$ in ${\frak M}_3$ contains a small cover as its
representative.
\end{thm}

 \begin{proof}
If $\beta=0$, then the bounding small cover $M^{'}$ in
Lemma~\ref{l8} can be chosen as a representative of $\beta$. If
$\beta\not=0$, then  $\beta$  is a linear combination of those 13
small covers stated in Theorem~\ref{dim}. Consider the small cover
$M$ constructed  in Lemma~\ref{l8} and take a fixed point $p$ of $M$
with representation $\rho_1\rho_2\rho_3$, one first preforms an
equivariant connected sum $M\sharp M$ of two copies of $M$ along the
fixed point $p$, so that $M\sharp M$ is also a small cover and
bounds equivariantly. Obviously, one can still preform an
equivariant connected sum of $M\sharp M$ with some chosen
arbitrarily from 13 small covers stated in Theorem~\ref{dim} such
that the resulting 2-torus manifold is a small cover.  This means
that $\beta$ must contain a small cover as its representative.
 \end{proof}

\end{document}